\documentclass[10pt]{article}
\usepackage{amsmath,amssymb,amsfonts,amsthm,mathtools,xfrac,extarrows,enumitem, bbold}
\usepackage{graphicx, float}
\usepackage{color}  
\usepackage[T1]{fontenc}  
\usepackage{amscd}  
\usepackage[active]{srcltx}    
\usepackage[font={small}]{caption} 
\usepackage{epstopdf}
\usepackage[T1]{fontenc}
\usepackage{amscd} 
\usepackage{xr,xr-hyper}

\numberwithin{equation}{section}  
\newcommand{\ov}{\overline}
\renewcommand{\a }{\alpha }
\renewcommand{\b }{\beta }
\renewcommand{\d}{\delta }
\renewcommand{\l}{\lambda }

\newcommand{\R}{\mathbb{R}}

\makeatletter
\renewcommand*\env@matrix[1][*\c@MaxMatrixCols c]{%
  \hskip -\arraycolsep
  \let\@ifnextchar\new@ifnextchar
  \array{#1}}
\makeatother

\newtheorem{corollary}{Corollary}[section]
\newtheorem{remark}{Remark}[section]    
\newtheorem{lemma}{Lemma}[section]
\newtheorem{proposition}{Proposition}[section] 
\newtheorem{theorem}{Theorem}
\newtheorem{definition}{Definition}[section]

\numberwithin{equation}{section}

\def\hline
{
\vspace{-12pt}
\begin{center}
\rule{0.7\linewidth}{0.5pt}
\end{center}
\vspace{12pt}
}

\newcommand{\intbar}{\mathop{\int\makebox(-13.5,0){\rule[4pt]{.7em}{0.3pt}}%
\kern-6pt}\nolimits}

\newcommand{\pa}{\partial}

\newcommand{\e }{\varepsilon }
\newcommand{\n }{\nabla }
\newcommand{\s }{\sigma }
\newcommand{\D }{\Delta }

\renewcommand{\a }{\alpha }
\renewcommand{\b }{\beta }
\renewcommand{\d}{\delta }
\renewcommand{\l}{\lambda }

\title{Prescribing Morse scalar curvatures: pinching and Morse theory}

\author
{
Andrea Malchiodi
\;\&\;
Martin Mayer
}
\begin{document} 

\maketitle

\begin{abstract} 
  
\noindent We consider the problem of prescribing conformally the scalar curvature 
on compact manifolds of positive Yamabe class in 
dimension $n \geq 5$. We prove new existence results 
using Morse theory and some analysis on blowing-up solutions, 
under suitable pinching conditions on the curvature function. We also provide  new non-existence results
showing the sharpness of some of our assumptions, both in terms of the dimension and 
of the Morse structure of the prescribed function. 
\end{abstract}

\tableofcontents

\section{Introduction}

We deal here with the classical problem of prescribing the scalar curvature of closed manifolds, 
whose study initiated systematically with the papers \cite{kw1}, \cite{kw2}, \cite{kw}. We will consider in particular {\em conformal} changes of metric.
On $(M^n,g_0),\,n \geq 3$ and for a smooth positive function $u$ on $M$ we denote by 
$$g = g_u = u^{\frac{4}{n-2}} g_{0}$$
a metric $g$ conformal to $g_0$. 
Then the scalar curvature transforms according to 
\begin{equation}
\label{eq:conf-change-curv}
R_{g_u} u^{\frac{n+2}{n-2}} 
= 
L_{g_0} u 
:=
- c_{n} \Delta_{g_{0}} u + R_{g_{0}} u,
\quad 
c_n = \frac{4(n-1)}{(n-2)}, 
\end{equation} 
see \cite{aul}, Chapter 5, \textsection 1, where $\Delta_{g_{0}}$ is the Laplace-Beltrami operator of $g_0$. The elliptic operator $L_{g_0}$ is known as the {\em conformal Laplacian} and obeys the covariance law  
\begin{equation}\label{eq:covar}
L_{g_{u}} (\phi) 
= 
u^{-\frac{n+2}{n-2}} L_{g_0}(u \phi)
\quad \text{ for } \quad
\phi \in C^{\infty}(M). 
\end{equation}
If under a conformal change of metric one wishes to prescribe the scalar curvature of $M$ as a given function $K : M \longrightarrow \R$, 
by \eqref{eq:conf-change-curv} one would then need to find positive solutions of the nonlinear elliptic problem 
\begin{equation}\label{eq:scin}
L_{g_0} u 
=
K u^{\frac{n+2}{n-2}}  
\quad \text{ on } \quad  
(M,g_0). 
\end{equation}
The above equation is variational and of {\em critical type}, and  it  presents 
a lack of compactness. When $K$ is zero or negative,  in which case $(M,g_0)$ has to 
be of zero or negative {\em Yamabe class} respectively, the nonlinear term in the equation 
makes the Euler-Lagrange energy for \eqref{eq:scin}  coercive and 
solutions always exist, as proved in \cite{kw} via the method of sub- and super solutions. 
In the same paper though Kazdan and Warner showed that for $K$ positive 
there are obstructions to existence. Indeed, if $f : S^n \longrightarrow \R$ is the 
restriction to the sphere of a coordinate function in $\R^{n+1}$, then 
\begin{equation}\label{eq:kw}
\int_{S^n}  \langle \nabla K, \nabla f \rangle_{g_{S^n}} u^{\frac{2n}{n-2}} \, d\mu_{g_{S^{n}}} 
= 
0, 
\end{equation}
for all solutions $u$ to \eqref{eq:scin}.  This forbids for example the prescription of 
affine functions or generally of functions $K$ on $S^n$ that are monotone 
in one Euclidean direction. More examples are given in \cite{be}.

Existence of solutions for $K$ positive on manifolds of positive Yamabe class were found some years later. 
In the spirit of a result by Moser in \cite{mo}, where antipodally symmetric curvatures were prescribed on $S^2$, 
in \cite{es} the authors showed solvability of \eqref{eq:scin} on $S^n$, when $K$ is 
invariant under a group of isometries without fixed points and satisfies suitable {\em flatness} 
assumptions depending on the dimension. Other results with symmetries were also found in 
\cite{hv}, \cite{hv93}. 

Another theorem, regarding more general functions $K$, was proved in \cite{bab} and \cite{bc} for the case of $S^3$ assuming that $K : S^3 \longrightarrow \R_{+}$ is a Morse function satisfying the 
generic condition  
\begin{equation}\label{eq:nd}
\{ \nabla K  = 0 \} \cap \{ \Delta K = 0 \} = \emptyset
\end{equation}
together with the {\em index formula} 
\begin{equation}\label{eq:bcin}
\sum_{\{ x \in M\; : \; \nabla K(x) = 0, \Delta K(x) < 0\}}(-1)^{m(K,x)} \neq (-1)^n, 
\end{equation}
where  $m(K,x)$ denotes the Morse index of $K$ at $x$, cf. 
\cite{[ACGPY]}, \cite{[ACPY1]}, \cite{[ACPY2]}, \cite{[SZ]}.
 
To put our work into context, it is useful to briefly describe the strategy to prove the 
latter result. A useful tool for studying \eqref{eq:scin} in the spirit of
\cite{su} is its \emph{subcritical approximation} 
\begin{equation}\label{eq:scin-tau}
L_{g_0} u  
=
K  u^{\frac{n+2}{n-2}-\tau},\; \quad
0<\tau \ll 1, 
\end{equation}
which up to  rescaling $u$ is the Euler-Lagrange equation for the functional  
\begin{equation}\label{eq:JKt}
J_{\tau}(u)
=
\frac
{\int_{M} \left( c_{n}\vert \nabla u\vert_{g_{0}}^{2}+R_{g_{0}}u^{2} \right) d\mu_{g_{0}}}
{(\int_M Ku^{p+1}d\mu_{g_{0}})^{\frac{2}{p+1}}},
\quad 
p = \frac{n+2}{n-2} - \tau. 
\end{equation}
By its scaling-invariance and the sign-preservation  of its gradient flow, we  assume $J_\tau$ to be defined on 
\begin{equation}\label{eq:XXX}
X 
=
\{
u\in W^{1,2}(M,g_{0})\mid u\geq 0\;\wedge\; \|u\| = 1
\}, 
\end{equation}
where the norm $\| \cdot \|$ is defined by \eqref{eq:norm} in case of a positive Yamabe class. 
The advantage of \eqref{eq:scin-tau} is that with a sub-critical exponent the problem is now compact and solutions can be easily found. On the other hand one might expect solutions to {\em blow-up} as $\tau \longrightarrow 0$. However, as for the above mentioned result, sometimes it is possible to 
completely classify blowing-up solutions and to show by degree- or Morse-theoretical arguments, that  
there must be solutions to \eqref{eq:scin-tau}, 
which do not blow-up 
and hence converge to solutions of \eqref{eq:scin}.  

When  blow-up occurs, there is a formation of {\em bubbles}, namely 
profiles that after a suitable dilation solve \eqref{eq:scin} on $S^n$ with $K\equiv 1$, cf. \cite{auw}, \cite{cgs}, \cite{tal}. 
In three dimensions due to a slow  decay, which implies that  mutual interactions among bubbles are {\em stronger} 
than the interactions of each bubble with $K$, it is possible to show that 
only one bubble can form at a time. Such bubbles develop necessarily at critical points of $K$ with negative 
Laplacian and their total contribution to the Leray-Schauder degree of \eqref{eq:scin-tau} is precisely the 
summand in \eqref{eq:bcin}, just taken with the opposite sign. Then by compactness of the equation and the Poincar\'e-Hopf theorem  the total degree of \eqref{eq:scin-tau} is 1, 
contradicting  inequality \eqref{eq:bcin}. In \cite{yy}, \cite{yy2} this result was extended to $S^n$ under suitable flatness conditions 
on $K$, which are similar to those in \cite{es}, cf. \cite{yy2}, \cite{[BCCH4]} for $K$ Morse with a formula 
different from \eqref{eq:bcin} on $S^4$, where only finitely-many blow-ups may occur, but only at restricted locations. 
Results of different kind were also proven in \cite{[CD]} for $n = 2$
and in \cite{[BI]}, \cite{[BCCHhigh]}, \cite{[BIEG]}, cf. Chapter 6 in \cite{aul}. 

\medskip

In higher dimensions the analysis of blowing-up solutions to \eqref{eq:scin-tau} for $\tau \longrightarrow 0$ 
is more difficult. Some results are available in \cite{cl1}-\cite{cl4}, showing that in general blow-ups 
with infinite energy may occur. For $K$ Morse on $S^n$ and still satisfying \eqref{eq:nd} and \eqref{eq:bcin} 
some results in general dimensions were proven under suitable 
{\em pinching conditions}, cf. \cite{agp}, \cite{ab}, \cite{[ACPY3]}, \cite{cxy}, \cite{cx12} and \cite{mal}.

\medskip 

In our first theorem we extend the result in \cite{cx12} to Einstein manifolds of 
positive Yamabe class under the  pinching condition
\begin{equation}\label{eq:p1} \tag{$P_1$}
\frac{K_{\max}}{K_{\min}} 
\leq  
2^{\frac{1}{n-2}}, 
\end{equation}
where with obvious notation 
$$
K_{\max} = \max_{S^n} K
\quad \text{ and } \quad
K_{\min} = \min_{S^n} K.
$$  
If $K$ is Morse, it must have a non-degenerate maximum and hence
\eqref{eq:bcin} requires the 
existence of at least a second critical point of $K$ with negative Laplacian. We also show that the existence of two such critical points is sufficient for existence under a more stringent pinching requirement, 
namely 
\begin{equation}\label{eq:p2} \tag{$P_2$}
\frac{K_{\max}}{K_{\min}} \leq  \left( \frac{3}{2} \right)^{\frac{1}{n-2}}. 
\end{equation}

\begin{theorem} \label{t:pinching}
Suppose $(M^n,g_0)$ is an Einstein manifold of positive Yamabe class with $n \geq 5$, and that $K$ is a positive 
Morse function on $M$ verifying \eqref{eq:nd}. Assume we are in one of the following two situations: 
\begin{enumerate}
\item[(i)] $K$ satisfies \eqref{eq:p1} and \eqref{eq:bcin}; 
\item[(ii)] $K$ satisfies \eqref{eq:p2} and has at least two critical points with negative Laplacian. 
\end{enumerate}
Then \eqref{eq:scin} has a positive solution. 
\footnote
{
In the case of $S^n$ the curvature pinching assumptions of 
Theorem 1 (i) are stronger than those of Thoerem 1.2 in \cite{cx12},  
but we cannot completely follow the proof there.  
We refer in particular to the continuity of $T_1$ before formula (7.6) in \cite{cx12}. Its definition depends on the quantity 
	$\| v - 1\|$, which tends to zero for every initial datum $u_0$ as an evolution time $t$ tends to infinity. 
	However, since the quantity $\| v - 1\|$ may not be globally monotone in time,  we are unable to verify the continuity of $T_1$. 
}
\end{theorem}

\medskip 

The pinching conditions we require can indeed be relaxed, even though they become 
more technical to state, see Theorem \ref{t:pinching-ref} for details.

\medskip

\begin{remark}\label{'r:intro} 
\begin{enumerate}[label=(\roman*)]
\item 
We would like to emphasize \cite{bab1} as the first work to 
analyse with a high degree of generality the lack of compactness of the conformally prescribed Morse scalar curvature problem on higher dimensional spheres and the first one to provide non trivial existence results, 
which are based on a  \textit{topological invariant} introduced by A. Bahri in the same work. This invariant might prove  useful in relaxing or even removing the pinching assumptions in Theorem 1.
\item 
Also in higher dimensions, but considering only the zero weak limit scenario,
we also refer to our previous work \cite{MM2} and \cite{M3} in the subcritical and critical case respectively for a comprehensive discussion of the aforementioned lack of compactness.
\item 	
To our knowledge condition (ii) is of new type and the restriction on the dimension is optimal. Building on some 
non-existence result in \cite{str05} for the Nirenberg problem on $S^2$, it is 
possible to manufacture curvature functions on $S^3$ and on $S^4$ such that
under condition  (ii), even under arbitrary pinching problem \eqref{eq:scin} has no solution, cf. Remark \ref{r:Struwe34}. 

Such  curvatures can be obtained  perturbing affine functions, forbidden by the 
Kazdan-Warner obstruction, and  deforming their non-degenerate maximum into two 
nearby maxima and a saddle point. In low dimension candidate solutions are ruled out 
via blow-up analysis, as they could form at most one bubble. A contradiction to  existence 
is then obtained by a {\em quantitative} version of \eqref{eq:kw}, showing that even if the 
integrand changes sign, the total integral does not vanish. In dimension $n \geq 5$ the 
contradiction argument breaks down, since multi-bubbling occurs, as 
shown in \cite{lz} for $n=6,7,8,9$, cf. \cite{cny}. 
\end{enumerate}

\end{remark}

We are going to describe next our strategy for proving Theorem \ref{t:pinching}, which 
relies on the  subcritical approximation \eqref{eq:scin-tau}. We considered in \cite{MM1} 
a special class of solutions to the latter equation, namely solutions with uniformly bounded energy and zero weak limit. 
Even though in high dimension general blow-ups, as described before, can have a 
complicated behaviour, we proved that this class of solutions can only develop {\em isolated simple 
ones}, i.e. at most one bubble per blow-up point, cf. Subsection \ref{ss:iiisbu} for precise definitions. These occur at critical points of $K$ with negative 
Laplacian with no further restriction on their location, as shown in \cite{MM2}, see also \cite{M3} and 
\cite{M5} for the relation with a dynamic approach to \eqref{eq:scin}. 

The outcome of these results, summarized in Theorem \ref{t:ex-multi}, is that if \eqref{eq:scin} is not solvable 
and  $(u_{\tau_n})_n$ is a sequence of solutions to \eqref{eq:scin-tau}  
with uniformly bounded energy as $\tau_n \longrightarrow 0$, then they are in one-to-one correspondence 
with the finite sets 
$$
\{ x_1, \dots, x_q \} \subseteq \{ \nabla K  = 0 \} \cap \{ \Delta K < 0 \}
,
\quad q \geq 1. 
$$
Such solutions $u_{\tau, x_1, \dots, x_q}$ are also non-degenerate  for the functional $J_\tau$ on $X$, cf. \eqref{eq:JKt}, \eqref{eq:XXX}, 
and their Morse index and asymptotic energy can be explicitly computed, depending on $(K(x_i))_i$ 
and on $(m(K,x_i))_i$. This allows then 
to deduce existence results via variational or Morse-theoretical arguments.

The stronger  the pinching  of $K$ is, the more the above solutions $u_{\tau, x_1, \dots, x_q}$ tend 
to {\em quantize} in energy, depending on the number of blow-up points. Energy sublevels 
of $J_\tau$ within these strata can then be deformed to sublevels of the {\em reference} subcritical 
Yamabe energy  
$\bar{J}_{\tau}$ defined on $X$ as 
\begin{equation*}
\bar{J}_{\tau}(u)
= 
\frac
{\int_{M} \left( c_{n}\vert \nabla u\vert_{g_{0}}^{2}+R_{g_{0}}u^{2} \right) d\mu_{g_{0}}}
{(\int_M u^{p+1}d\mu_{g_{0}})^{\frac{2}{p+1}}}. 
\end{equation*}
It turns out that on Einstein manifolds the only critical points of $\bar{J}_{\tau}$ are constant functions, cf. Theorem 6.1 in 
\cite{bvv}, and therefore all sublevels of $\bar{J}_{\tau}$ are contractible. The pinching 
condition allows to show that suitable sublevels of $J_\tau$ are also contractible. As a consequence 
 the total degree of single-bubbling solutions is equal to one, while 
the total degree of doubly-bubbling solutions, which must occur at couples of distinct points in 
$\{ \n K = 0\} \cap \{ \Delta K < 0 \}$,
is equal to zero. By direct computation we can then  deduce existence of solutions under both conditions (i) and (ii) 
in Theorem \ref{t:pinching}.

\medskip

One may wonder whether  stronger pinching  assumptions might  induce existence 
under weaker conditions than the second one in (ii). In view of the Kazdan-Warner obstruction and 
of Remark \ref{'r:intro},  it is tempting to think that when $n \geq 5$ and  $K : S^n \longrightarrow \R_{+}$ has 
more than just one local maximum and minimum, solutions may always exist. We show that 
this is not the case, and that critical points of $K$ with positive Laplacian are less relevant. 
For $K$ Morse on $S^n$ we define  
\begin{equation}\label{eq:MjK}
\mathcal{M}_j(K) 
=
\sharp \left\{ x \in S^{n} \: : \; \nabla K (x) = 0\; \wedge \; m(K,x) = j  \right\}. 
\end{equation}
We then have the following result.

\begin{theorem} \label{t:non-ex}
For $n \geq 3$ and any Morse function $\tilde{K} : S^n \longrightarrow \R_{+}$  with only one local maximum point, there exists 
a Morse function $K : S^n \longrightarrow \R$ such that 
\begin{enumerate}[label=(\roman*)]
\item  $\mathcal{M}_j(K) = \mathcal{M}_j(\tilde{K})$ for all $j$;
\item the Laplacian at all critical points of $K$ with the exception of its 
local maximum is positive;
\item there is no conformal metric on $S^n$ with scalar curvature 
$K$.
\end{enumerate}
$K$ can be also chosen so that $\frac{{K}_{\max}}{{K}_{\min}}$ 
is arbitrarily close to $1$. 
\end{theorem}

\begin{remark} In comparison to the latter result we note, that the non-existence examples in \cite{be} for $S^2$ 
	are not pinched and imply the existence of one or more local maxima. 
	\end{remark}

Theorem \ref{t:non-ex} is proved by composing curvature functions 
as those discussed in Remark \ref{'r:intro} (iii) with a reflection with respect to the last Euclidean coordinate.  We construct a suitable sequence of curvatures ${K}_m$ 
as in Theorem \ref{t:non-ex} converging to a monotone function in the last Euclidean variable of $\R^{n+1} \supseteq S^n$
with a non-degenerate maximum at the north pole and all other critical points, with positive Laplacian, 
accumulating near the south pole of $S^n$.

Assuming by contradiction that \eqref{eq:scin} 
has solutions $u_m$ with $K = {K}_m$, by a result in \cite{cl1}, \cite{cl97} such solutions would stay uniformly 
bounded away from both poles. As we noticed before,  blow-ups in high 
dimensions might have diverging energy. However, near the south pole both the mutual interactions 
among bubbles and that of each bubble with ${K}_m$ would tend to \emph{deconcentrate} 
highly-peaked solutions. Via some Pohozaev type identities, this can be made rigorous showing 
first that blow-ups at the south pole are \emph{isolated simple} and then that they indeed do not occur. 
The delicate part in this step is that the critical point structure of $(K_m)_m$ is degenerating, and 
we still need uniform controls on solutions.

The analysis near the north pole is harder, since  the two  interactions just described have competing effects. 
We need then to rule out different limiting scenarios for sequences of candidate solutions, namely regular limits, 
singular limits and zero  limits locally away from the north pole. The latter case is the most delicate: we show that a 
regular bubble must form at  a slowest possible blow-up rate and via Kelvin inversions, decay estimates and 
integral identities, that blow-up cannot occur.

\

Our strategy also allows to improve some existing results in the literature with assumptions 
that are  {\em localized} in the range of $K$, as for example in \cite{[BI]}, cf.\cite{[ACPY1]}, \cite{[ACPY2]} 
and \cite{str05}
for $n=2$. The general idea is to use min-max schemes, e.g. the mountain pass, and to 
 use competing paths whose maximal energy lies below that of every possible blowing-up solution 
for \eqref{eq:scin-tau} with bounded energy, via the pinching conditions. The fact that such blow-ups are isolated simple 
reduces the number of diverging competitors,  permitting us to relax  previous pinching constraints in the literature. We can 
also use Morse-theoretical arguments, in particular {\em relative Morse inequalities}, to prove existence  by {\em counting} the number 
of min-max paths and  of diverging competitors, cf.  Subsection \ref{ss:p-mm}.

	\
	
The plan of the paper is the following: in Section \ref{s:prel} we collect some 
preliminary material on the variational structure of the problem, on singular 
solutions to the Yamabe equation and on blow-up analysis. In Section \ref{s:pinch} 
we prove  existence results via index counting or min-max theory, exploiting the 
pinching conditions. In 
Section \ref{ss:n-e} we then prove  non-existence results by constructing suitable curvature functions with 
prescribed Morse structure and using blow-up analysis to find contradiction to existence. We finally collect the proofs of some technical 
results in an appendix.

\

\section{Preliminaries} \label{s:prel}

In this section we gather some background and preliminary material concerning the 
variational structure of the problem, with a description of subcritical bubbling with finite energy. 
We also collect some integral identities, the notion of simple blow-up and some of its consequences,
as well as  some properties of singular Yamabe metrics.

\subsection{Variational structure}\label{ss:var}

We consider a  closed Riemannian manifold  
$
M=(M^{n},g_{0}) 
$
with induced volume measure $\mu_{g_{0}}$ and scalar curvature $R_{g_{0}}$.
For $X$  as in \eqref{eq:XXX}  
the {\em Yamabe invariant}  is 
\begin{equation*}\begin{split}
Y(M,g_{0})
= &
\inf_{u \in X}
\frac
{\int \left( c_{n}\vert \nabla u \vert_{g_{0}}^{2}+R_{g_{0}}u^{2} \right) d\mu_{g_{0}}}
{(\int u^{\frac{2n}{n-2}}d\mu_{g_{0}})^{\frac{n-2}{n}}},
\quad 
c_{n}=4\frac{n-1}{n-2}, 
\end{split}\end{equation*}
which due to \eqref{eq:conf-change-curv} 
depends only on the conformal class of $g_0$. 
We will restrict ourselves to manifolds of \emph{positive Yamabe class}, namely those 
for which the Yamabe invariant is positive. 
In this case the {\em conformal Laplacian} 
$ 
L_{g_{0}}=-c_{n}\Delta_{g_{0}}+R_{g_{0}}
$ 
is a positive and self-adjoint operator and admits a Green's function 
$$
G_{g_{0}}:M\times M \setminus \Delta \longrightarrow\R_{+},
$$
where $\D$ is the diagonal of $M \times M$. 
For a conformal metric 
$$g=g_{u}=u^{\frac{4}{n-2}}g_{0}$$ 
there holds
\begin{equation*}\begin{split} d\mu_{g_{u}}=u^{\frac{2n}{n-2}}d\mu_{g_{0}}
\quad \text{ and } \quad 
R=R_{g_{u}}
=
u^{-\frac{n+2}{n-2}}(-c_{n} \Delta_{g_{0}} u+R_{g_{0}}u) 
=
u^{-\frac{n+2}{n-2}}L_{g_{0}} u, 
\end{split}\end{equation*}
and by the positivity of $L_{g_0}$ there exist constants $c, C > 0$ such that 
\begin{align*}
c\Vert u \Vert_{W^{1,2}(M,g_0)}^2 
\leq
\int u \, L_{g_{0}}u \, d\mu_{g_{0}}
=
\int \left( c_{n}\vert \nabla u \vert^{2}_{g_{0}}+R_{g_{0}}u^{2}
\right) d\mu_{g_{0}}
\leq 
C\Vert u \Vert_{W^{1,2}(M,g_0)}^2. 
\end{align*}
Therefore the square root of 
\begin{equation}\label{eq:norm}
\Vert u \Vert^2 
=
\Vert u \Vert_{L_{g_0}}^2 = \int u \, L_{g_{0}}u \, d\mu_{g_{0}}
\end{equation}
can be used  as an equivalent norm on $W^{1,2}(M,g_0)$. Setting  
$$
R=R_{u} \quad \text{ for } \quad g=g_{u}=u^{\frac{4}{n-2}}g_{0}
$$ 
we have 
\begin{equation}\label{eq:r}
r=r_{u}
=
\int R d\mu_{g_{u}}
=
\int u L_{g_{0}} ud\mu_{g_{0}} 
\end{equation}
and hence from \eqref{eq:JKt}
\begin{equation}\label{eq:kp}
J_{\tau}(u)
=
\frac{r}{k_{\tau}^{\frac{2}{p+1}}}
\quad \text{ with } \quad 
k_{\tau} 
=
\int K  u^{p+1} d \mu_{g_{0}}. 
\end{equation}
The first- and second-order derivatives of the functional $J_\tau$ are given by 
\begin{equation}\label{first_variation_evaluating} 
\partial J_{\tau}(u)v
= 
\frac{2}{k_{\tau}^{\frac{2}{p+1}}}
\big[\int L_{g_{0}}uvd\mu_{g_{0}}-\frac{r}{k_{\tau}}\int Ku^{p}vd\mu_{g_{0}}\big], 
\end{equation}  
and 
\begin{equation}\label{second_variation_evaluating} 
\begin{split}
\partial^{2} J_{\tau}(u) vw
= &
\frac{2}{k_{\tau}^{\frac{2}{p+1}}}
\big[\int L_{g_{0}}vwd\mu_{g_{0}}-p\frac{r}{k_{\tau}}\int Ku^{p-1}vwd\mu_{g_{0}}\big] \\
& -
\frac{4}{k_{\tau}^{\frac{2}{p+1}+1}}
\big[
\int L_{g_{0}}uvd\mu_{g_{0}}\int Ku^{p}wd\mu_{g_{0}} \\
& \quad\quad\quad\quad\quad+
\int L_{g_{0}}uwd\mu_{g_{0}}\int Ku^{p}vd\mu_{g_{0}}
\big] \\
& +
\frac{2(p+3)r}{k_{\tau}^{\frac{2}{p+1}+2}}
\int Ku^{p}vd\mu_{g_{0}}\int Ku^{p}wd\mu_{g_{0}}.
\end{split}
\end{equation}
Note that $J_\tau$ is scaling-invariant in $u$, whence we may restrict our attention to $X$, see \eqref{eq:XXX}. 
$J_{\tau}$ is of class $C^{2, \alpha}_{\text{loc}}$ and its critical points, suitably scaled, give rise to 
solutions of \eqref{eq:scin-tau}. Furthermore its $L_{g_{0}}$- gradient flow 
preserves the condition $\vert \cdot \vert = 1$ as well as non-negativity of initial data, in particular the set $X$.

\subsection{Finite-energy bubbling}\label{ss:feb}

\noindent
{\em Bubbles} denote concentrated solutions  of \eqref{eq:scin} or \eqref{eq:scin-tau} with the profile of conformal factors of Yamabe metrics 
on $S^n$. We follow our notation from \cite{MM1}, \cite{MM2}.

Let us recall  the construction of {\em conformal normal coordinates} from \cite{lp}. 
Given $a \in M$, these are  geodesic normal coordinates 
for a suitable conformal metric $g_{a} \in [g_{0}]$. 
If $r_a$ is  the geodesic distance from $a$ with respect to the metric $g_a$, 
the expansion of the Green's function  for the conformal 
Laplacian $L_{g_{a}}$ with pole at $a \in M$, denoted by 
$G_{a}=G_{g_{a}}(a,\cdot)$, simplifies considerably. From Section 6 of \cite{lp}
\begin{equation}\begin{split} \label{eq:exp-G}
G_{ a }
=
\frac{1}{4n(n-1)\omega _{n}}(r^{2-n}_{a}+H_{ a }), \quad r_{a}=d_{g_{a}}(a, \cdot)
, \quad 
H_{ a }=H_{r,a }+H_{s, a }
\end{split}\end{equation}
for $g_{a}=u_{a}^{\frac{4}{n-2}}g_{0}$.
Here $H_{r,a }\in C^{2, \alpha}_{\text{loc}}$ is a {\em regular} part, while the {\em singular one}  is of type 
\begin{equation*}
\begin{split}
H_{s,a}
=
O
\begin{pmatrix}
r_{a} & \text{ for }\quad n=5
\\
\ln r_{a} & \text{ for }\quad n=6
\\
r_{a}^{6-n} & \text{ for } \quad n\geq 7
\end{pmatrix}.
\end{split}
\end{equation*} 
For $\l > 0$ large let us define
\begin{equation}\label{eq:bubbles}
\varphi_{a, \lambda }
= 
u_{ a }\left(\frac{\lambda}{1+\lambda^{2} \gamma_{n}G^{\frac{2}{2-n}}_{ a }}\right)^{\frac{n-2}{2}}, \quad 
G_{ a }
=
G_{g_{ a }}( a, \cdot), \quad
\gamma_{n}=(4n(n-1)\omega_{n})^{\frac{2}{2-n}}. 
\end{equation}
The constant $\gamma_{n}$ is chosen in order to have 
$$
\gamma_{n}G^{\frac{2}{2-n}}_{ a }(x) 
= 
d_{g_a}^2(a,x) + o(d_{g_a}^2(a,x))
\quad \text{ as } \quad
x \longrightarrow a.
$$
Rescaled by a suitable factor depending on $K(a)$, for large values of $\l$ the functions $\varphi_{a, \lambda }$ are approximate solutions of \eqref{eq:scin}; 
moreover for $\l^{-2} \simeq \tau$ they are also approximate solutions to \eqref{eq:scin-tau} since in this 
regime $\l^{-\tau} \longrightarrow 1$ as $\tau \longrightarrow 0$, cf. Theorem \ref{t:ex-multi} below. Up a scaling constant their profile 
is given by the function
\begin{equation}\label{eq:U0}
U_0(x) = (1+\vert x\vert ^{2})^{\frac{2-n}{2}}
\quad \text{ for } \quad 
x \in \R^n,
\end{equation}
cf. Section 5 in \cite{MM1},
which realizes the best constant in the Sobolev inequality,  i.e.
\begin{equation}\label{eq:Sn}
\hat{c}_0 
=
c_n \inf_{0 \not\equiv u \in C_c^\infty(\R^n)} 
\frac
{\int_{\R^n} \vert \n u\vert ^2 dx}
{\left( \int_{\R^n} \vert u\vert ^{2^*} dx \right)^{\frac{2}{2^*}}} 
=
c_n 
\frac
{\left({\Gamma (n)}/{\Gamma \left(\frac{n}{2}\right)}\right)^{2/n}}
{\pi  (n-2) n},
\quad 
2^* = \frac{2n}{n-2}.
\end{equation}

\medskip 
\noindent {\bf{Notation.}} 
For a finite set of points $(x_i)_i$ in $M$ and $$K  : M \longrightarrow \R$$ 
a Morse function we will use the short notation 
\begin{equation}\label{eq:not}
K_i = K(x_i) \quad \text{ and } \quad m_i = m(K,x_i). 
\end{equation}
Combining the main results in \cite{MM1} and \cite{MM2} one has the following theorem. 
\begin{theorem}\label{t:ex-multi} (\cite{MM1}, \cite{MM2})
Let $(M,g)$ be a closed manifold of dimension $n \geq 5$ of positive Yamabe class and 
$K : M \longrightarrow \R$ be a positive Morse function satisfying \eqref{eq:nd}. 
Let $x_1, \dots, x_q$ be distinct critical points of $K$ with negative Laplacian. Then, as $\tau \longrightarrow 0$, 
there exists a unique solution $u_{\tau, x_1, \dots, x_q}$ developing exactly one bubble 
at each point $x_i$ and converging weakly to zero in $W^{1,2}(M,g)$ as $\tau \longrightarrow 0$. 
\\
Precisely there exist $\l_{1,\tau}, \dots, \l_{q,\tau} \simeq \tau^{-\frac 12}$ and points 
$a_{i,\tau} \longrightarrow x_i$ for all $i$ 
such that 
$$
\left\| u_{\tau, x_1, \dots, x_q} - \sum_{i=1}^q K_i^{\frac{2-n}{4}}  \varphi_{a_{i,\tau}, \lambda_{i,\tau}} \right\| \longrightarrow 0
\quad \text{ and } \quad
J_\tau(u_{\tau, x_1, \dots, x_q} ) 
\longrightarrow 
\hat{c}_0 \left( \sum_{i=1}^q  K_i^{\frac{2-n}{2}} 
\right)^{\frac{2}{n}} 
$$
as $\tau \longrightarrow 0$. Up to scaling $u_{\tau, x_1, \dots, x_q}$ is non-degenerate for $J_\tau$ 
and  
$$
m(J_\tau, u_{\tau, x_1, \dots, x_q}) 
= 
(q-1) + \sum_{i=1}^q 
(n-m_i)
.
$$ 
Conversely all blow-ups of \eqref{eq:scin-tau} with uniformly bounded energy and zero weak limit 
are as above.   
\end{theorem}
In \cite{MM1}, \cite{MM2} we proved  much more precise asymptotics on the 
solutions provided above, which are not  needed here, 
but were useful to show non-degeneracy. 
Recall also that the above statement 
is false for $n \leq 4$ since in three dimensions there could be at most one blow-up (in fact, 
no blow-up at all if $(M,g_0)$ is not conformally equivalent to $(S^3, g_{S^3})$ by the results in 
\cite{lizhu}), while 
in four dimensions there are constraints on blow-up configurations depending on $K$ and 
on the Green's function of $L_{g_0}$, cf. \cite{[BCCH4]} and \cite{yy2}.

\subsection{Integral identities and isolated simple blow-ups}\label{ss:iiisbu}

For finite-energy blow-ups of \eqref{eq:scin} one can prove a decomposition of solutions 
into finitely-many bubbles  in the spirit of \cite{str84}, see Section 3 in \cite{MM1}. 
In Section \ref{ss:n-e}  we will deal instead with general solutions, and some 
tools and definitions will be useful in this respect. 

\medskip 

\noindent
Recall Pohozaev's identity in a Euclidean ball $B_r = B_r(0) \subseteq \R^n$
for solutions to 
\begin{equation}\label{eq:K}
- c_n \D  u 
= 
K u^{\frac{n+2}{n-2}} 
\quad \text{ in } \quad 
\ov B_r. 
\end{equation}
If $\nu$ is the outer unit normal to $\pa B_r$, solutions of this equation satisfy 
\begin{equation}\label{eq:poho}
\frac{1}{2^*}  \int_{B_r} \sum_i x_i \frac{\pa K}{\pa x_i} u^{2^*} dx 
= 
\frac{1}{2^*} \oint_{\pa B_r} \langle x, \nu \rangle K u^{2^*} d \s 
+ 
c_n \oint_{\pa B_r} B(r, x, u, \n u) \, d \s, 
\end{equation}
where 
\begin{equation}\label{eq:BBB}
B(r, x, u, \n u) 
= 
\frac{n-2}{2} u \frac{\pa u}{\pa \nu} - \frac{1}{2}  \langle x, \nu \rangle
\vert \n u\vert ^2 +\frac{\pa u}{\pa \nu} \langle \nabla u, x \rangle. 
\end{equation}
This well-known identity is derived multiplying the equation by $x_i \frac{\pa u}{\pa x_i}$ and 
integrating by parts, cf. Corollary 1.1 in \cite{yy}. We describe next a 
{\em translational version} of it. Multiply \eqref{eq:K} by $\frac{\pa u}{\pa x_i}$ 
to get 
$$
- c_n \int_{B_r} \frac{\pa u}{\pa x_i} \D u \, dx 
= 
\frac{1}{2^*} \int_{B_r} K (u^{2^*})_{x_i} dx. 
$$
By the Gauss-Green theorem this becomes 
\begin{equation*}
\begin{split}
-  c_n\oint_{\pa B_r} \frac{\pa u}{\pa x_j} \langle \nu, e_j \rangle \frac{\pa u}{\pa x_i} d \s 
& + 
\frac{1}{2} c_n \int_{B_r} (\vert \n u\vert ^2)_{x_i} dx  \\
= & 
\frac{1}{2^*} \oint_{\pa B_r} 
K u^{2^*} \langle \nu, e_i \rangle d \s 
- 
\frac{1}{2^*} \int_{B_r} u^{2^*} \frac{\pa K}{\pa x_i} dx, 
\end{split}
\end{equation*}
where $e_j$ denotes the $j$-th standard basis vector of $\R^n$.
\begin{lemma}\label{l:poho-transl}
Let $u$ solve \eqref{eq:K} in $\ov B_r$ with $K \in C^1(\ov B_r)$. 
Then for all $i = 1, \dots, n$ 
\begin{equation}\label{eq:poho-trans}
\begin{split}
- c_n \oint_{\pa B_r} \frac{\pa u}{\pa x_j} \langle \nu, e_j \rangle \frac{\pa u}{\pa x_i} d \s 
& + 
\frac{1}{2} c_n \oint_{\pa B_r} \vert \n u\vert ^2 \langle \nu, e_i \rangle d \s \\
= &
\frac{1}{2^*} \oint_{\pa B_r} K u^{2^*} \langle \nu, e_i \rangle d \s 
- 
\frac{1}{2^*} \int_{B_r} u^{2^*} \frac{\pa K}{\pa x_i} dx. 
\end{split}
\end{equation}

\end{lemma}

\medskip

Consider now a sequence $(u_m)_m$ of solutions to 
\begin{equation}\label{eq:Ki}
- c_n \D  u_m 
= 
K_m(x) u_m^{\frac{n+2}{n-2}} 
\quad \text{ in } \quad 
\ov B_r, 
\quad \text{ with } \quad 
u_m(x_m) \longrightarrow \infty.  
\end{equation}
If $x_m \longrightarrow \bar x \in M$, the point 
$\bar x$ is called a {\em blow-up point} for $(u_m)_m$. 
For $r > 0$ let
$$
\overline{u}_m(r) = \intbar_{\pa B_r(x_m)} u_m d \s
$$
denote the radial average and we define 
\begin{equation}\label{eq:barwm}
\overline{w}_m(r) = r^{\frac{n-2}{2}} \overline{u}_m(r). 
\end{equation}
Following standard terminology, we define convenient classes of blow-ups.

\begin{definition}\label{d:isol-simple} Let $\xi_m$ be a local maximum for $u_m$. A blow-up point $\bar \xi = \lim_m \xi_m$
for $u_m$ is said to be \emph{isolated} if there exist (fixed)  constants $\rho > 0$ and $C > 0$ 
such that for all $m$ large 
\begin{equation}\label{eq:isol}
u_m(x) \leq \frac{C}{\vert x - \xi_m\vert ^{\frac{n-2}{2}}} 
\quad \text{ for } \quad
\vert x - \xi_m\vert  \leq \rho. 
\end{equation}
The blow-up point is said to be \emph{isolated  simple} if there exists $\rho \in (0,\infty]$ (fixed) 
such that for all $m$ large $\overline{w}_m(r)$ has precisely one critical point in $(0,\rho)$. 
\end{definition}

\medskip 

\noindent 
The above definitions are useful to characterize  {\em bubble towers} and single bubbles respectively,
yielding convergence after dilation and further estimates.  
If $(K_m)_m$ is a sequence of 
positive functions uniformly bounded in $C^1(\ov B_r)$ and bounded away from zero, 
we have the next result on isolated simple blow-ups, which is a consequence of Proposition 
2.3 in \cite{yy}.

\begin{lemma} \label{l:sing-prof} 
Suppose that $u_m$ solves \eqref{eq:Ki} with 
$$
C^{-1} \leq K_m \leq C 
\quad \text{ and } \quad 
\vert \n K_m\vert  \leq C,\; C > 0$$ 
and that $0 \in B_r$ is  
an isolated simple blow-up. Then there exists $C > 0$ such that 
\begin{equation}\label{eq:ub-um}
u_m(x) \leq C \, u_m(\xi_m)^{-1} \vert x - \xi_m\vert^{2-n} 
\quad \text{ in } \quad
B_{r/2}(\xi_m). 
\end{equation}
Moreover in a fixed neighbourhood $U$ of zero one has  
$$
u_m(\xi_m) u_m(x) 
\longrightarrow 
z(x) = a \vert x\vert ^{2-n} + h(x)
\quad \text{ in } \quad
C^2_{\text{loc}}(U \setminus \{0\}), 
$$
where 
$z>0$ 
is singular harmonic on $U \setminus \{0\}$,
$a > 0$ constant and $h$ smooth and harmonic at $x=0$. 
\end{lemma}

We first remark that after a suitable blow-down procedure $U$ can possibly coincide with all of $\R^n$, in which case $h$ has to be identically constant and non-negative. Secondly, the same holds true if $U$ coincides with $\R^n$ minus a discrete 
set $S$ of points including the origin and 
$$
z(x) 
= 
\sum_{p_i \in S} a_i \vert x-p_i\vert ^{2-n} + \tilde{h}(x),
$$ 
in which case $\tilde{h}$ is  constant.  
We can then apply  (2)  of Proposition 1.1 in \cite{yy} to conclude that for $r > 0$ small, 
if $0$ is an isolated simple blow-up, then 
\begin{equation}\label{eq:lim-B}
\oint_{\pa B_r} B(r, x, u_m, \n u_m) d \s 
= 
- \frac{(n-2)^2}{2}  \frac{h(0) \omega_n}{r^{n-2} u_m(\xi_m)^{2}} 
(1 + o_m(1) + o_r (1)), 
\end{equation}
where $\omega_n = \vert S^{n-1}\vert $, $h$ is as in Lemma \ref{l:sing-prof} and  $o_m(1) \xlongrightarrow{m\to \infty} 0$, $o_r(1) \xrightarrow{r\to 0} 0$. 

\medskip 

We next recall the following well-known result which can be found in \cite{snotes} and stated  in  Section 8 of \cite{kms}. It follows 
by iteratively extracting bubbles from  solutions large in  $L^\infty$-norm. 

\begin{proposition} \label{p:kms}
Consider on $S^n$  a function $K : S^n \longrightarrow \R_{+}$  satisfying for 
$$
C_0^{-1} \leq K \leq C_0
\quad \text{ and } \quad 
\|K\|_{C^2(S^n)} \leq C_0
$$
some $C_0 \geq 1$. 
Given $\d > 0$ small and $R > 0$ large, there exists 
$C = C(\d,R,C_0) > 0$ 
such that, if $u$  
solves \eqref{eq:scin} with such $K$ and  $\max_{S^n} u \geq C$, then there exist local maxima
$\xi_1, \dots, \xi_N \in S^n$, $N = N(u) \geq 1$ of $u$ such that 
\begin{enumerate}[label=(\roman*)]
\item \quad the balls $(B_{r_i}(\xi_i))_{i=1}^N$ with $r_i = R u(\xi_i)^{-\frac{2}{n-2}}$ are disjoint; 
\item \quad  in normal coordinates $x$ at $\xi_i$ one has 
$$
\left\|
u(\xi_i)^{-1} u(u(\xi_i)^{-\frac{2}{n-2}} y) 
- 
\left( 1 + k_i \vert  y \vert ^2 \right)^{\frac{2-n}{2}}
\right\|_{C^2(B_R(0))} 
< 
\d, 
$$ 

\quad
where $k_i = \frac{1}{n(n-2)c_n} K(\xi_i)$ and  $y = u(\xi_i)^{\frac{2}{n-2}} x$; 
\item \quad 
$
u(x) 
\leq 
C d_{S^n}(x,\{\xi_1, \dots, \xi_N\})^{-\frac{n-2}{2}}
$ 
for all $x \in S^2$;
\item \quad
$
d_{S^n}(\xi_i,\xi_j)^{\frac{n-2}{2}} u(\xi_j) 
\geq 
C^{-1}
$ 
for all $i \neq j$.
\end{enumerate}
\end{proposition}

\subsection{Singular solutions and conservation laws}\label{ss:sing-sol}

We recall next some properties of radial \emph{singular solutions} (at $x = 0$) of the critical equation 
$$
- \Delta u 
= 
\kappa \, u^{\frac{n+2}{n-2}} 
\quad \text{ in } \quad
\R^n \setminus \{0\}
\quad \text{ with } \quad 
\kappa > 0. 
$$
Such solutions are of interest as they could arise as limits of regular solutions, see 
Theorem 1.4 in \cite{cl2}. By 
Theorem 8.1 in \cite{cgs} all the singular solutions of the above equation are radial,  
cf. \cite{mp} for other properties. 
If we  look for solutions in the form 
$$
u(x) = \vert x\vert ^{\frac{2-n}{2}} v(\log \vert x\vert ), 
$$
then by direct computation $v$  satisfies
$$
- v''(t) 
= 
\kappa  v^{\frac{n+2}{n-2}}(t) 
- 
\left( \frac{n-2}{2} \right)^2 v (t). 
$$
The latter is a Newton equation of the form $v''(t) = - V'(v(t))$, with potential 
$$
V(v) 
= 
\kappa  \frac{n-2}{2n} v^{\frac{2n}{n-2}} 
- 
\frac{1}{2}\left( \frac{n-2}{2} \right)^2  v^2. 
$$
This implies the conservation of the  {\em Hamiltonian energy} 
$$
\frac{1}{2} (v')^2 
+ 
\kappa\frac{n-2}{2n} v^{\frac{2n}{n-2}} 
- 
\frac{1}{2} \left( \frac{n-2}{2} \right)^2  v^2  
=: 
H.  
$$
The value 
\begin{equation*}
v_{0} 
\equiv 
\left[\left( \frac{n-2}{2} \right)^2 \kappa^{-1} \right]^{\frac{n-2}{4}}
\quad \text{ with  Hamiltonian} \quad 
H_{0}
=
- \frac{1}{n} \kappa \left[\left( \frac{n-2}{2} \right)^2 \kappa^{-1} \right]^{\frac{n}{2}}
\end{equation*}
is the only critical point of $V$ on the 
positive $v$-axis and for every value 
$
H 
\in 
\left( 
H_{0}
, 
0 
\right)
$ 
there is a unique positive periodic solution $v_H$,   
called {\em Fowler's solution}, with period increasing in $H$ and tending 
to infinity as $H \longrightarrow 0$. 
In fact, as $H \longrightarrow 0$, $v_H$ converges on the compact sets of $\R$ to a 
homoclinic solution $v_0$ tending to zero for $t \longrightarrow \pm \infty$, where 
 $v_0$ corresponds to a {\em regular solution} to the above Yamabe equation.

\begin{lemma}\label{l:cons-law}
For 
$
H 
\in 
\left( 
H_{0}
, 
0 
\right)
$
let $u_H(x) = \vert x\vert ^{\frac{2-n}{2}} v_H(\log \vert x\vert )$.
Then 
$$
\frac{1}{2^*} \oint_{\pa B_r} \langle x, \nu \rangle \kappa \, u_H^{2^*} d \s 
+ 
\oint_{\pa B_r} B(r, x, u_H, \n u_H) \, d \s 
= 
\omega_n H, 
$$ 
where $\omega_n = \vert S^{n-1}\vert $ and $B$ is as in \eqref{eq:BBB}. 
\end{lemma}

\begin{proof} 
In terms of $u_H, u'_H$, after some cancellation the boundary integrand becomes
$$
\frac{1}{2} r (u'_H)^2 
+ 
\frac{1}{n} (n - 2)/2 \kappa r u_H^{\frac{2n}{n-2}} 
+ 
\frac{n-2}{2}  u_H u'_H. 
$$
We have clearly that  
$$
\n u_H(x) 
= 
\frac{2-n}{2 \vert x\vert }  u_H(\vert x\vert ) 
+ 
\vert x\vert ^{\frac{2-n}{2}-1} v'_H(\log \vert x\vert ). 
$$
Substituting for $v_H$, the boundary integrand transforms into  
$$
r^{1-n} \left(\left(4 n (v'_H)^2-(n-2)^2 n v_H^2\right)
+ 
4 (n-2)\kappa v_H^{\frac{2 n}{n-2}}\right)
= 
8nr^{1-n} H. 
$$
Integrating on $\pa B_r$, the conclusion immediately follows. 
\end{proof}

\section{Existence results}\label{s:pinch}

In this section we prove Theorem \ref{t:pinching} and other existence results, using pinching assumptions on $K$ and Morse-theoretical arguments.

\subsection{Pinching and topology of sublevels} \label{ss:pinch-top}
 
Here we show that a suitable pinching condition implies contractibility in $X$ of some sublevels of $J_\tau$
for $(M,g_0)$ Einstein. Such conditions will be made more explicit in the next subsection, depending on the critical points of $K$.  Recall that $p = \frac{n+2}{n-2} - \tau$ and $K:M\longrightarrow \R_{+}$ is strictly positive and let 
\begin{equation*}
\begin{split}
\overline{A}= \left( \frac{K_{\max}}{K_{\min}} \right)^{\frac{2}{p+1}}  \underline{A}
\quad \text{ for any } \quad
\underline{A} > 0. 
\end{split}
\end{equation*}

\begin{proposition}\label{p:gap}
Let $(M^n,g_0)$ be an Einstein manifold of positive Yamabe class and $\tau > 0$. 
If
\begin{equation*}
\begin{split}
\{\partial J_{\tau}=0\} 
\;\cap\;
\{
\underline A \leq J_{\tau} \leq \overline A
\}
=
\emptyset
\end{split}
\end{equation*}
for some $\underline{A}>0$, then for every
$c \in [\underline A, \overline A]$ 
the sublevel $\{J_\tau \leq c\}$ is contractible. 
\end{proposition}

\begin{proof}  For $u \in X$ we  clearly have 
$$
K_{\max}^{-\frac{2}{p+1}} \bar{J}_\tau (u) \leq J_\tau(u) 
\leq 
K_{\min}^{-\frac{2}{p+1}} \bar{J}_\tau (u), 
$$
whence for $A, B > 0$
$$
J_\tau(u) \leq A  \Longrightarrow  \bar{J}_\tau(u) \leq K_{\max}^{\frac{2}{p+1}} 
A
\quad \text{ and } \quad 
\bar{J}_\tau(u) \leq B  \Longrightarrow  J_\tau(u) \leq K_{\min}^{-\frac{2}{p+1}} B. 
$$ 
Therefore we have the for $\underline{A} > 0$ inclusions
$$
\{ J_\tau \leq \underline A \} 
\subseteq 
\{ \bar{J}_\tau \leq K_{\max}^{\frac{2}{p+1}} \underline{A} \}
\subseteq 
\{ J_\tau \leq \overline A  \}.  
$$
As
$\partial J_\tau$ is uniformly bounded on sublevels and of class $C^{1,\a}$ there, cf. \eqref{first_variation_evaluating}, \eqref{second_variation_evaluating}, 
the negative gradient flow $\phi$ for $J_\tau$ 
with respect to the scalar product induced by $L_{g_{0}}$ 
is globally well defined on $X$, see \eqref{eq:XXX}, and in time and $\phi(t,u)$ depends continuously on the initial condition $u$. 
Note that $\phi$  preserves the $L_{g_{0}}$-norm, see \eqref{eq:norm}, as well as non-negativity of initial 
data and hence the set $X$, cf. Section 4 in \cite{M3}. 

\

\begin{figure}[h]\label{fig_1}
\begin{center}
\includegraphics[width=0.8\textwidth]{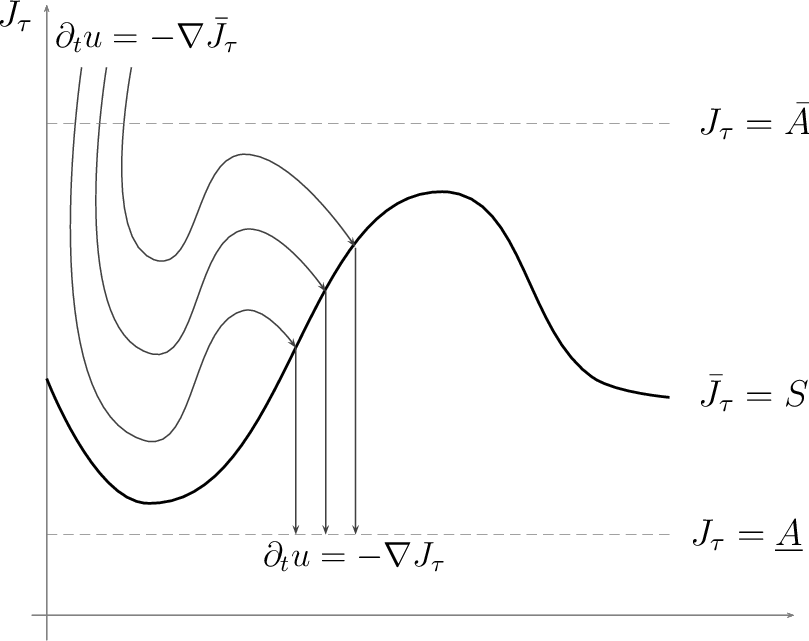}
\caption{$(\tau)$-Yamabe and prescribed scalar curvature flows combined}
\end{center}
\end{figure}

\noindent  
Since
$$
\{\partial J_{\tau}=0\} 
\;\cap\;
\{
\underline A \leq J_{\tau} \leq \overline A
\}
=
\emptyset
$$
and $J_\tau$ satisfies the Palais-Smale condition, as $\tau > 0$, by the deformation lemma, cf. Section 7.4 in \cite{am} and \emph{transversality} 
for any 
$u 
\in 
[\underline A \leq J_{\tau} \leq \overline A]$ 
there exists a first time $T_u \geq 0$, which is continuous in $u$, such that  
$$
\phi(T_u, u) \in \left\{ J_\tau \leq \underline{A} \right\}. 
$$
Recalling that 
$\{ \bar{J}_\tau \leq K_{\max}^{\frac{2}{p+1}} \underline{A} \} 
\subseteq 
\{ J_\tau \leq \overline A  \}$, 
consider then the homotopy
$$
F 
: 
[0,1] \times (\{ \bar{J}_\tau \leq K_{\max}^{\frac{2}{p+1}} \underline{A} \}) 
\longrightarrow 
X
:  
(s,u) \longmapsto \phi(s \, T_u, u). 
$$ 
If $u$ belongs to the sublevel $\{ J_\tau \leq \underline{A} \}$, 
then $T_u = 0$ and hence 
$$F(s,u) = u \quad \text{  for all  } \quad s \in [0,1].$$ 

\

\noindent
Therefore $F$ deforms
$\{ \bar{J}_\tau \leq K_{\max}^{\frac{2}{p+1}} \underline{A} \}$ 
into 
$\{ J_\tau \leq \underline A \}$, 
but not necessarily within 
$$\{ \bar{J}_\tau \leq K_{\max}^{\frac{2}{p+1}} \underline{A} \}.$$ 
This can be achieved composing $\phi$ on the left with a suitable {\em Yamabe-type flow}. 
Recall  that, if $(M^n,g_0)$ is Einstein and of positive Yamabe class, by Theorem 6.1 in \cite{bvv} the equation  
$$
L_{g_0} u = u^{p}, \quad p = \frac{n+2}{n-2} - \tau
$$ 
has only constant solutions. Hence the  
infimum of $\bar{J}_\tau$ is attained and equal to 
$R_{g_0} Vol_{g_0}(M)^{1-\frac{2}{p+1}}$. 
Since  the Palais-Smale condition holds also for $\bar{J}_\tau$, the gradient 
flow $\bar \phi(t,u)$ of $\bar{J}_\tau$ evolves all initial data $u$ to a 
constant function, intersecting \emph{transversally} every level set of $\bar{J}_\tau$ higher 
than its infimum. Similarly to the previous reasoning there  
exists for any $u \in X$ a first time $\overline{T}_u \geq 0$, continuous in $u$, such that  
$$
\overline{\phi}(\overline{T}_u, u) 
\in 
\left\{ \bar{J}_\tau \leq K_{\max}^{\frac{2}{p+1}} \underline{A} \right\}. 
$$
Defining 
$$
\tilde{F}(s,u) 
= 
\overline{\phi}(\overline{T}_{F(s,u)}, F(s,u)),  
$$
we deduce that $\tilde F(s,u)$ is a \emph{deformation retract} of 
$\{ \bar{J}_\tau \leq K_{\max}^{\frac{2}{p+1}} \underline{A} \}$ 
onto 
$\{ J_\tau \leq \underline A \}$
and therefore realizes a homotopy equivalence, cf. Chapter II in \cite{mas}.

\

\noindent
On the other hand every non-empty sublevel 
$\{ \bar{J}_\tau \leq B \}$,
in particular 
$\{ \bar{J}_\tau \leq K_{\max}^{\frac{2}{p+1}} \underline{A} \}$, 
is via the deformation lemma and Palais-Smale's condition for $\bar{J}_\tau$  homotopically equivalent to a point. 
Hence we deduce the same property for
$\{ J_\tau \leq \underline A \}$.
Still by the deformation lemma and the Palais-Smale condition,  this is true also for 
$\{ J_\tau \leq c \}$ 
with 
$c \in [\underline A, \overline A]$.
This concludes the proof. 
\end{proof}

\medskip

\noindent 
For the above proof to work, it is indeed sufficient to assume that the functional $J_{\tau}$ for $\tau=0$  
has no critical points in a restricted energy range.

\subsection{Pinching and degree counting} \label{ss:pdc}
 
If problem \eqref{eq:scin} has no solutions, using Theorem \ref{t:ex-multi}  
we will show that Proposition \ref{p:gap} applies, provided suitable 
pinching conditions on $K$ hold true.  Arguing by contradiction, we will then derive existence results of which Theorem \ref{t:pinching} is a particular case. To that end we first order the set
$$
\{ x_1, \dots, x_l \} 
= 
\{ \n K = 0 \} \cap \{ \D K < 0 \}
$$ 
so that 
\begin{equation*}
K_{1}=K(x_{1})\geq \ldots \geq K_{l}=K(x_{l}). 
\end{equation*}
Recalling our notation in \eqref{eq:Sn} and  \eqref{eq:not}, for $m \in \{1, \dots, l\}$ we then define  
\begin{equation}\label{eq:unovEm}
\underline{E}_m 
= 
\hat{c}_0 \left( \sum_{i=1}^m  K_i^{\frac{2-n}{2}} 
\right)^{\frac{2}{n}} 
\quad \text{ and } \quad 
\overline{E}_m 
=  
\hat{c}_0
\left( \sum_{i=l-m+1}^l  K_i^{\frac{2-n}{2}} 
\right)^{\frac{2}{n}}. 
\end{equation}
As we will see, these numbers represent the minimal and maximal limit energies for solutions 
developing $m$ bubbles and weakly converging to zero as $\tau \longrightarrow 0$. 
We then have the following result.

\begin{proposition} \label{p:pinching=>gap}
Suppose that \eqref{eq:scin} has no solutions, and assume that 
\begin{equation}\label{eq:tPm} \tag{$\tilde P_m$}
\left( \frac{K_{\max}}{K_{\min}} \right)^{\frac{n-2}{2}} 
< 
\frac{\underline{E}_{m+1}}{\overline{E}_m} 
\end{equation}
 for some 
$m \in \{1, \dots, l-1\}$. Then there exists $0<\e \ll 1$ such that 
\begin{equation*}
\Big\{ \partial J_{\tau}=0 \Big\}
\; \cap \;
\Big\{
(1 + \e) \overline{E}_m
\leq 
J_{\tau} 
\leq
\left( \frac{K_{\max}}{K_{\min}} \right)^{\frac{2}{p+1}} (1 + \e) \overline{E}_m  
\Big\}
=
\emptyset
\end{equation*}
for all $\tau>0$ sufficiently small.
\end{proposition}
 \begin{proof} 
Suppose \eqref{eq:scin} has no positive solutions. Then, as $\tau \searrow 0$, all positive solutions of 
\eqref{eq:scin-tau} with 
uniformly bounded energy must have zero weak limit. These are then described by Theorem  \ref{t:ex-multi} and  of the form 
$u_{\tau, x_{i_1}, \dots, x_{i_q}}$ 
with $x_{i_1}, \dots, x_{i_q}$ distinct points of $\{ x_1, \dots, x_l \}$
and energies
\begin{equation*}
J_{\tau}(u_{\tau, x_{i_1}, \dots, x_{i_q}})
\longrightarrow
\hat{c}_0 \left( \sum_{j=1}^q  K_{i_j}^{\frac{2-n}{2}} \right)^{\frac{2}{n}}
\quad \text{ as } \quad
\tau \longrightarrow 0. 
\end{equation*}
By the way we ordered the points $(x_i)_i$, we clearly have that 
$$
\hat{c}_0 \left( \sum_{j=1}^q  K_{i_j}^{\frac{2-n}{2}}\right)^{\frac{2}{n}} 
\leq 
\overline{E}_m 
\quad \text{ for } \quad
q \leq m
$$
and
$$
\hat{c}_0 \left( \sum_{j=1}^q  K_{i_j}^{\frac{2-n}{2}} 
\right)^{\frac{2}{n}} 
\geq 
\underline{E}_{m+1} 
\quad \text{ for } \quad
q \geq m+1 
$$ 
Then the statement immediately follows.  
\end{proof}

\begin{remark}\label{r:rel-pinching}
Let us consider the  pinching condition 
\begin{equation} \label{eq:Pm} \tag{$P_m$}
\frac{K_{\max}}{K_{\min}}
<
\left(\frac{m+1}{m}\right)^{\frac{1}{n-2}}.
\end{equation}
We then have
\begin{equation}\label{eq:Pm-to-tPm}
(P_{m+1}) \; \Longrightarrow \;  (P_{m}) 
\quad \text{ and } \quad 
(P_{m}) \; \Longrightarrow \;  (\tilde P_{m}) 
\quad \text{ for all } \quad 
m \geq 1.  
\end{equation}
Indeed, while the first implication is obvious, for the second we find from \eqref{eq:Pm} 
\begin{equation*}
\sum_{i=1}^{m+1}K_{i}^{\frac{2-n}{2}}
\geq 
\frac{m+1}{K_{\max}^{\frac{n-2}{2}}}
>
(\frac{K_{\max}}{K_{\min}})^{\frac{n-2}{2}}\frac{m}{K_{\min}^{\frac{n-2}{2}}}
\geq 
(\frac{K_{\max}}{K_{\min}})^{\frac{n-2}{2}}
\sum_{i=l-m+1}^{l}K_{i}^{\frac{2-n}{2}}, 
\end{equation*}
which implies  \eqref{eq:tPm} by the definitions in \eqref{eq:unovEm}. 
Finally we observe that also 
\begin{equation}\label{eq:tPm1m2}
(\tilde P_{m_1}) 
\; \Longrightarrow \;  
(\tilde P_{m_2}) 
\quad \text{ for all } \quad
m_1 \geq m_2.  
\end{equation}
Indeed we may argue inductively and see that $(\tilde P_{m_1})$ for 
$m_{1}=m_{2}+1$ 
implies 
\begin{equation*}
\begin{split}
\sum_{i=1}^{m_{2}+1}K_{i}^{\frac{2-n}{2}}
= &
\sum_{i=1}^{m_{1}+1}K_{i}^{\frac{2-n}{2}}
-
K_{m_{1}+1}^{\frac{2-n}{2}}  
>
(\frac{K_{\max}}{K_{\min}})^{\frac{n-2}{2}}
\sum_{i=l-m_{1}+1}^{l}K_{i}^{\frac{2-n}{2}}
-
K_{m_{1}+1}^{\frac{2-n}{2}}
\\
= &
(\frac{K_{\max}}{K_{\min}})^{\frac{n-2}{2}}
\sum_{i=l-m_{2}+1}^{l}K_{i}^{\frac{2-n}{2}}
+
(\frac{K_{\max}}{K_{\min}})^{\frac{n-2}{2}}K_{l-m_{1}+1}^{\frac{2-n}{2}}
-
K_{m_{1}+1}^{\frac{2-n}{2}}, 
\end{split}
\end{equation*}
and
\begin{equation*}
(\frac{K_{\max}}{K_{\min}})^{\frac{n-2}{2}}K_{l-m_{1}+1}^{\frac{2-n}{2}}
-
K_{m_{1}+1}^{\frac{2-n}{2}}
=
K_{l-m_{1}+1}^{\frac{2-n}{2}}
\left(
(\frac{K_{\max}}{K_{\min}})^{\frac{n-2}{2}}
-
(\frac{K_{l-m_{1}+1}}{K_{m_{1}+1}})^{\frac{n-2}{2}}
\right)
\geq 0.
\end{equation*}
We therefore obtain \eqref{eq:tPm1m2} as desired. 
\end{remark}
 
\medskip
 
We prove next the following result, which by \eqref{eq:Pm-to-tPm} in the previous remark extends 
Theorem \ref{t:pinching}. 

\begin{theorem} \label{t:pinching-ref}
Suppose $(M^n,g_0)$ is an Einstein manifold of positive Yamabe class with $n \geq 5$ and  $K$ is a positive 
Morse function on $M$ satisfying \eqref{eq:nd}. Assume we are in one of the following two situations: 
\begin{enumerate}
\item[(j)]  
$K$ satisfies $(\tilde{P}_1)$ and \eqref{eq:bcin}; 
\item[(jj)] 
$K$ satisfies $(\tilde{P}_2)$ and has at least two critical points with negative Laplacian. 
\end{enumerate}
Then \eqref{eq:scin} has a positive solution.
\end{theorem}

\begin{proof} The proof will be carried out by
contradiction, assuming that the functional $J_0$ does not have any critical 
point, so we have the conclusion of Proposition \ref{p:pinching=>gap}
and thus the conclusion of Proposition \ref{p:gap}.

	Suppose (j) holds: recalling \eqref{eq:unovEm}, we deduce that 
for $\e > 0$ small the sublevel 
$\{ J_\tau \leq (1 + \e) \overline{E}_1 \}$
is contractible and that $J_\tau$ has no critical points at level  
$(1 + \e) \overline{E}_1$. 
By Theorem \ref{t:ex-multi}, all critical points of $J_\tau$ at 
lower levels are single-bubbling solutions $u_{\tau,x_i}$, 
which totally contribute to the Leray-Schauder degree of \eqref{eq:scin-tau}  
by the amount 
$$
\sum_{ x_{j}\in \{\nabla K=0\}\cap \{\Delta K<0\}} (-1)^{n-m_j},
$$
see  \eqref{eq:not}. 
By the Poincar\'e-Hopf theorem this total sum must be equal to the Euler characteristic 
$\chi(\{ J_\tau \leq (1 + \e) \overline{E}_1 \}) = 1$, 
which contradicts the assumption. 

\medskip

Suppose now that (jj) holds true, and let us again assume  that $J_0$ has no critical points. As 
$(\tilde{P}_2)$ implies $(\tilde{P}_1)$, see Remark \ref{r:rel-pinching}, 
we thus have a contradiction from case (j), provided \eqref{eq:bcin} holds. 
Hence we may assume that $(\tilde{P}_2)$  holds and 
\begin{equation}\label{eq:no-bcin}
\sum_{x_i \in \{ \nabla K = 0 \} \cap \{ \Delta K < 0\}}(-1)^{m_i}
= 
(-1)^n. 
\end{equation}
With the same reasoning as above we obtain that for $\e > 0$ small the sublevel 
$\{ J_\tau \leq (1 + \e) \overline{E}_2 \}$ 
is contractible and that $J_\tau$ has no critical points at level  
$(1 + \e) \overline{E}_2$.

By our assumptions solutions of \eqref{eq:scin-tau} with limiting energies less or equal to 
$(1 + \e) \overline{E}_2$ 
are either single- or doubly-bubbling solutions. By \eqref{eq:no-bcin} the contribution 
of the former to the degree is 1, while the contribution of the latter must be zero.

By Theorem \ref{t:ex-multi} doubly-bubbling solutions blow-up at 
\underline{distinct} critical points of $K$ with negative Laplacian, whence
by the characterization 
of their Morse index necessarily 
\begin{equation*}
0
=
\sum_{x_{i}\neq x_{j}, x_{i},x_{j}\in \{\nabla K=0\}\cap \{\Delta K<0\}}
(-1)^{n-m_i+n-m_j}. 
\end{equation*}
Combining the last formula with \eqref{eq:no-bcin}, we compute
\begin{equation*}
\begin{split}
0
= &
\sum_{x_{i}\neq x_{j}, x_{i},x_{j}\in \{\nabla K=0\}\cap \{\Delta K<0\}}
(-1)^{n-m_i+n-m_j} \\
= &
\sum_{x_{i}\in \{\nabla K=0\}\cap \{\Delta K<0\} }(-1)^{n-m_i}
\sum_{x_{j}\neq x_{i}, x_{j}\in \{\nabla K=0\}\cap \{\Delta K<0\}}
(-1)^{n-m_j}
\\
= &
\sum_{x_{i}\in \{\nabla K=0\}\cap \{\Delta K<0\} }(-1)^{n-m_i}
[
-(-1)^{n-m_i}
+
\sum_{ x_{j}\in \{\nabla K=0\}\cap \{\Delta K<0\}}(-1)^{n-m_j}
].
\end{split}
\end{equation*}
Using \eqref{eq:no-bcin} for the latter sum, we get 
\begin{equation*}
\begin{split}
0
= &
\sum_{x_{i}\in \{\nabla K=0\}\cap \{\Delta K<0\} }(-1)^{n-m_i}
[-(-1)^{n-m_i}
+
1 
] \\
= &
-
\sum_{x_{i}\in \{\nabla K=0\}\cap \{\Delta K<0\} }(-1)^{2(n-m_i)}
+
\sum_{x_{i}\in \{\nabla K=0\}\cap \{\Delta K<0\} }(-1)^{n-m_i}.
\end{split}
\end{equation*}
Again we know that the latter sum equals $1$, consequently
\begin{equation*}
\begin{split}
0
=  &
-\sum_{x_{i}\in \{\nabla K=0\}\cap \{\Delta K<0\} }(-1)^{2(n-m_i)}
+
1 
= 
-
\sharp\left(\{\nabla K=0\}\cap \{\Delta K<0\}\right) 
+
1, 
\end{split}
\end{equation*}
where $\sharp$ denotes the cardinality. Hence we reach a contradiction once more. 
\end{proof}

 \begin{remark}
\begin{enumerate}[label=\arabic*)]
\item 
The restriction on the dimension for condition (jj) is sharp, 
 cf. Remark \ref{r:Struwe34} for details.  Our proof indeed relies on Theorem \ref{t:ex-multi}, which only holds in dimension $n \geq 5$.
\item 
One could replace the degree-counting argument 
by  Morse's inequalities. This was done in \cite{[SZ]} in three dimensions and in \cite{cx12} 
 in arbitrary dimension under suitable pinching conditions.  
\item 
Formula \eqref{eq:bcin} arises from computing the contribution to the degree 
 of all single-bubbling solutions. 
 Considering  
 the blowing-up solutions in Theorem \ref{t:ex-multi} and the Morse-index 
 formula there, 
 it can be easily seen 
 that the total degree 
 of multi-bubbling solutions is $1$. If \eqref{eq:scin} is not solvable, Proposition \ref{p:gap} could then be applied for 
 large values of $A$, since $J_\tau$ would have only finitely-many 
 solutions with bounded energy, but we would derive no useful information from the Poincar\'e-Hopf theorem. 
\item 
Condition (j) (resp., (jj)) is used to find sublevels of  $\bar J_\tau$ that contain 
every blowing-up solution of \eqref{eq:scin-tau} forming one bubble (resp., two bubbles) but 
not containing any solution forming two (resp., three) bubbles or more.  Further pinching 
restrictions does not seem to lead to different existence results,  in view of Theorem \ref{t:non-ex}. 
\item 
The argument of the proof allows to also show that the solution provided by the above theorem 
is a critical point of $J_{\tau}$ for $\tau=0$ below a given energy value, see the comment after
Proposition \ref{p:gap}. This value can be any number exceeding 
the limiting energy of doubly-bubbling or triply-bubbling solutions as in Theorem \ref{t:ex-multi}. 
The existence result is also stable under small perturbations of the Einstein metric and might extend 
to conformal classes of metrics with a unique Yamabe representative, cf. \cite{dpz}. 
\end{enumerate}
\end{remark}

\subsection{Pinching and min-max theory}\label{ss:p-mm}

Here we show how Theorem \ref{t:ex-multi} can be used to 
improve results in the literature that rely on min-max theory, cf.
\cite{[CD]}, \cite{[ACPY2]} and \cite{str05} in two dimensions or  
\cite{[BI]}. 
Also with this approach and under some circumstances  the pinching assumption in Theorem \ref{t:pinching} 
can be relaxed. We have first the following general result, which will be later 
specialized to simpler situations or variants.  
 
\begin{theorem}\label{t:pq}
Let $(M^n,g)$, $n \geq 5$ be a closed Riemannian manifold of positive Yamabe class and $K$
be a positive Morse function on $M$  satisfying \eqref{eq:nd}. Assume that 
there is a set $\Xi \subseteq M$ with $\mathcal{C}$ components 
that contains $p$ local maxima $x_1, \dots, x_p$ of $K$ 
and such that 
$$
\max_\Xi \, (K^{\frac{2-n}{2}}) 
< 
\min 
\left\{ 
\left( 
K(x_i)^{\frac{2-n}{2}} + K(x_j)^{\frac{2-n}{2}} 
\right)^{\frac{2}{n}} 
\; : \; 
x_i \neq x_j \text{ local maxima of } K  
\right\}. 
$$ 
Assume also that $K$ has $q \geq 0$ critical points of index $1$ in the range 
$$[\min_\Xi K, \max_{i \in \{1, \dots, p\}} K(x_i)).$$ Then 
\eqref{eq:scin} has a solution provided that $q < p - \mathcal{C}$. 
\end{theorem}
 
  \medskip

\begin{remark}
Following our proof, the above result and thence the other ones in this subsection can be extended 
to $S^3$ without any pinching requirement due to single-bubbling. Note that from \cite{lizhu} problem \eqref{eq:scin} 
is always solvable on other three-manifolds. In four dimensions one can 
relax the pinching condition using constraints on multi-bubbling solutions as found in 
\cite{[BCCH4]} and \cite{yy2}. 
\end{remark}

Before proving  Theorem \ref{t:pq} we need some preliminaries. 
First we specify more precisely the asymptotic profile of the \underline{single-bubbling} solutions 
$u_{\tau,x_i}$ as in Theorem \ref{t:ex-multi}. If $\varphi_{a,\l}$ is as in \eqref{eq:bubbles}, 
then there exists 
$$
a_{i,\tau} \longrightarrow x_i\in \{\nabla K=0\} \cap \{\Delta K<0\}
\quad \text{ and } \quad 
\l_{i,\tau}^2 
= 
- (1 + o_\tau(1)) c_2 \frac{\Delta K(x_i)}{K(x_i) \tau}
$$
as $\tau\longrightarrow 0$, where
$c_2 = c_2(n) > 0$, 
see Section 3 in \cite{MM2}, 
such that 
\begin{equation}\label{eq:ui-close}
\left\| u_{\tau,x_i} - \varphi_{a_{i,\tau},\l_{i,\tau}} \right\| 
= 
o_\tau(1). 
\end{equation}
We then map $\Xi \subseteq M$ as in Theorem \ref{t:pq} into the variational space 
$X \subseteq W^{1,2}$, 
cf.\eqref{eq:XXX}, in such a way that each point $x_i$ is mapped to  $u_{\tau,x_i}$,
and derive an upper bound on $J_\tau$ under the image of this an embedding. 
Precisely consider for $r_0 > 0$ smooth 
$$
\tilde{\l} : M \longrightarrow \R_{+}
\quad \text{ and } \quad 
\tilde{a} : M \longrightarrow M
$$ 
satisfying with $a_{i,\tau}$ and $\l_{i,\tau}$ as in \eqref{eq:ui-close}
$$
\begin{cases}
\tilde{\l} = \tau^{-1/2} 
& 
\text{ in } \quad 
M \setminus \cup_{i=1}^p B_{4r_0}(x_i); 
\\ 
\tilde{\l} = \l_{i,\tau} 
& 
\text{ in } \quad
B_{2r_0}(x_i)
\\ 
c \tau^{-1/2} \leq \tilde{\l} \leq C \tau^{-1/2} 
& 
\text{ in } \quad
M
\end{cases}  
$$
and 
$$
\begin{cases}
\tilde{a}(x) = x & 
\text{ in } \quad
M \setminus \cup_{i=1}^p B_{4r_0}(x_i); 
\\ 
\tilde{a}(x) = a_{i,\tau} 
& 
\text{ in } \quad
B_{2r_0}(x_i); 
\\ 
\tilde{a} \in B_{4r_0}(x_i) 
& 
\text{ in } \quad
B_{4r_0}(x_i) 
\end{cases} 
$$
for some fixed constants $0 < c < C$.  Finally let for $x \in \Xi$
\begin{equation}\label{eq:var-test}
\tilde{\varphi}_{x,\tau} 
=
\begin{cases}
\varphi_{\tilde{a}(x), \tilde{\l}(x)}  
& 
\text{ in } \quad
M \setminus \cup_{i=1}^p B_{2r_0}(x_i); 
\\ 
( 1 - \frac{d(x,x_i)}{2 r_0} ) u_{\tau,x_i} + \frac{d(x,x_i)}{2 r_0} 
\varphi_{a_{i,\tau}, \l_{i,\tau}}  
& 
\text{ in } \quad
B_{2r_0}(x_i).
\end{cases} 
\end{equation}
We then have the following result. 
 
\begin{lemma}\label{l:en-est}
If $\tilde{\varphi}_{x,\tau}$ is  as in \eqref{eq:var-test} and if $\hat{c}_0$ is given in \eqref{eq:Sn}, one has that 
$$
\sup_{x \in \Xi} J_\tau(\tilde{\varphi}_{x,\tau}/\| \tilde{\varphi}_{x,\tau} \|)  
\leq  
\hat{c}_0 \max_\Xi \, (K^{\frac{2-n}{2}})  
+ 
o_\tau(1) + o_{r_0}(1), 
$$ 
where $o_\tau(1) \longrightarrow 0$ as $\tau \searrow 0$ and $o_{r_0}(1) \xrightarrow{r_{0}\to 0} 0$. 
\end{lemma}

\begin{proof}
Since $J_\tau$ is uniformly Lipschitz on finite energy sublevels and is 
scaling invariant, by \eqref{eq:ui-close} we are reduced to 
prove that 
$$
J_\tau({\varphi}_{x,\tilde{\l}(x)})  
\leq 
\hat{c}_0  K(x)^{\frac{2-n}{n}}  
+ 
o_\tau(1) 
\quad \text{ as } \quad 
\tau \searrow 0.  
$$
To show this, note that 
$\varphi_{x, \tilde{\l}(x)}$
is bounded from above and below by powers of 
$$\tilde{\l}(x) \simeq \tau^{-1/2},$$ 
and that 
$\tilde{\l}(x)^\tau \longrightarrow 1$ 
as 
$\tau \longrightarrow 0$,
whence
$$
J_\tau({\varphi}_{x,\tilde{\l}(x)})  
= 
\frac
{
\int_{M} \left( c_{n}\vert \nabla {\varphi}_{x,\tilde{\l}(x)}\vert_{g_{0}}^{2}+R_{g_{0}}{\varphi}_{x,\tilde{\l}(x)}^{2} \right) d\mu_{g_{0}}
}
{
(\int_M K {\varphi}_{x,\tilde{\l}(x)}^{2^*}d\mu_{g_{0}})^{\frac{2}{2^*}}
} 
+ 
o_\tau(1), 
\quad 
2^*= \frac{2n}{n-2}. 
$$
Using a change of variables, it is easy to see that
\begin{equation*}
\begin{split}
\frac
{\int_{M} \left( c_{n}\vert \nabla {\varphi}_{x,\tilde{\l}(x)}\vert_{g_{0}}^{2}+R_{g_{0}}{\varphi}_{x,\tilde{\l}(x)}^{2} \right) d\mu_{g_{0}}}
{(\int_M K {\varphi}_{x,\tilde{\l}(x)}^{2^*}d\mu_{g_{0}})^{\frac{2}{2^*}}}  
=  & 
c_n K(x)^{\frac{2-n}{n}} 
\frac
{\int_{\R^n} \vert \n U_0\vert ^2 dx}
{\left( \int_{\R^n} \vert U_0\vert ^{2^*} dx \right)^{\frac{2}{2^*}}} 
+ 
o_\tau(1) \\ 
= &
\hat{c}_0 K(x)^{\frac{2-n}{n}}  
+ 
o_\tau(1), 
\end{split}
\end{equation*}
where $U_0$ is given by \eqref{eq:U0}. This concludes the proof. 
\end{proof}
 
 \medskip
 
\begin{proof} {\em of Theorem \ref{t:pq}.}
Arguing by contradiction, assume that \eqref{eq:scin} has no solutions. Then, as  noticed in the 
previous subsection, all solutions of \eqref{eq:scin-tau} with uniformly bounded energy 
must have zero weak limit. Fix $\e > 0$ small:  
we know by  Theorem \ref{t:ex-multi} that $J_\tau$ has at least $p$ local 
minima of the form 
$u_{\tau,x_1}, \dots, u_{\tau,x_p}$ 
such that for $\tau$ small there holds 
$$ 
J_\tau(u_{\tau,x_j}) 
< 
\hat{c}_0 (\min_{i \in \{1, \dots, p\}} K(x_i) )^{\frac{2-n}{n}} 
+ 
\e
$$ 
and such that, for all sufficiently small values of $\tau$, 
$J_\tau$ has no critical point at level 
$$\hat{c}_0 (\min_{i \in \{1, \dots, p\}} K(x_i) )^{\frac{2-n}{n}} + \e.$$ 
We can assume that for  $\tau$  small there is no critical point of $ J_\tau$ at level 
$$\hat{c}_0 \max_\Xi \, (K^{\frac{2-n}{2}}) + \e$$
and we can modify $J_\tau$ near all its local minima at level less or equal to 
$$\hat{c}_0 \max_\Xi \, (K^{\frac{2-n}{2}}) + \e,$$ 
which are non-degenerate by Theorem \ref{t:ex-multi}, in order to still have the Palais-Smale condition, to not generate new critical points and so that the modified minima are at level zero. Call 
$\tilde J_\tau$ 
the resulting functional, which we can take of class $C^{2,\alpha}$ as the original one, see Figure 2. It will also possess at least $p$ critical points at level zero.

\begin{figure}[h]\label{fig_2}
\begin{center} 
\includegraphics[width=0.8\textwidth]{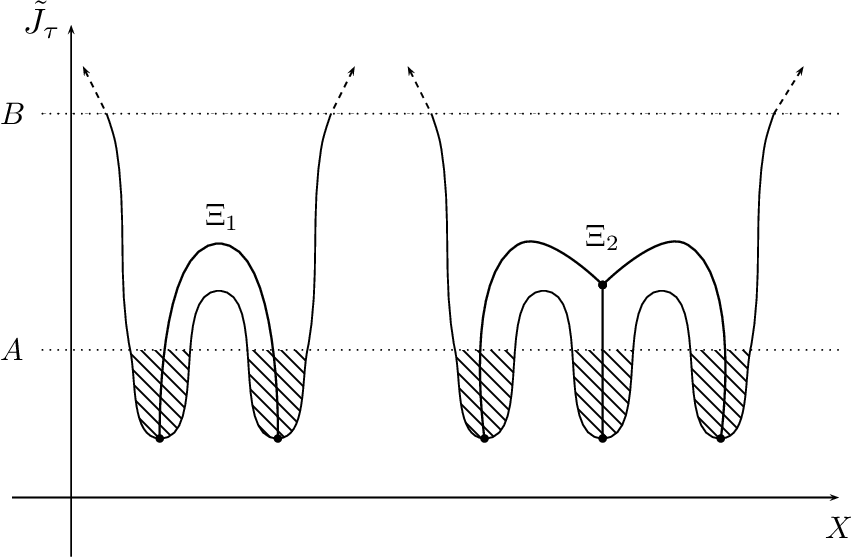}  
\caption{The modified functional $\tilde{J}_\tau$ and its sublevels.}
\end{center}
\end{figure}
We then use {\em relative Morse inequalities} for $\tilde{J}_\tau$, cf. Theorem 4.3 in \cite{Chang-book}, between the levels 
$$
A = \e
\quad \text{ and } \quad 
B = \hat{c}_0 \max_\Xi \, (K^{\frac{2-n}{2}}) + \e. 
$$
By construction $\tilde J_\tau$ has $C_0 = 0$ critical points of index zero 
and $C_1 = q$ critical points in the range $[A,B]$. 
Since $\tilde{J}_\tau$ has no local minima in the range $[A,B]$
and the Palais-Smale condition holds true, every point of
$\{ \tilde{J}_\tau \leq B \}$ 
can be joined to 
$\{ \tilde{J}_\tau \leq A \}$.
As a consequence  
$$
\beta_0 
:= 
\text{rank} \, H_0(\{ \tilde{J}_\tau \leq B \}, \{ \tilde{J}_\tau \leq A \}) 
= 
0, 
$$
see e.g. \cite{Hat}, Chapter 2, Exercise 16, page 130. On the other hand consider 
$$
\beta_1 
:= 
\text{rank} \, H_1(\{ \tilde{J}_\tau \leq B \}, \{ \tilde{J}_\tau \leq A \}). 
$$
Recall that 
$$
H_1(\{ \tilde{J}_\tau \leq B \}, \{ \tilde{J}_\tau \leq A \}) 
= 
Z_1(\{ \tilde{J}_\tau \leq B \}, \{ \tilde{J}_\tau \leq A \}) 
/ 
B_1(\{ \tilde{J}_\tau \leq B \}, \{ \tilde{J}_\tau \leq A \}), 
$$
where $Z_1$ and $B_1$ denote kernel and image of the boundary operators in one and 
two homological dimensions respectively,
cf. \cite{mas}, Chapter VII, \textsection 6. We claim next 
\begin{equation}\label{beta_1_inequality}
\beta_1 \geq p - \mathcal{C}.
\end{equation} 
To prove this, let 
$\Xi_1, \dots, \Xi_{\mathcal{C}}$ 
denote the connected components of $\Xi$. 
As our assumptions improve or stay invariant if we remove 
components containing none or only one point among $x_1, \dots, x_p$, 
we can assume that each component of $\Xi$ contains at least two among the points $x_1, \dots, x_p$. 

\

\noindent
Given $\Xi_j$ let 
$$X_j = \{ x_{i_1}, \dots x_{i_{\mathcal{C}_j}} \}$$ 
denote the local maxima of $K$ belonging to $\Xi_j$. Considering a curve 
$$
\gamma_{j,l} : [0,1] \longrightarrow M
\quad \text{ with } \quad 
\gamma(0)=x_{i_1}
\quad \text{ and } \quad  
\gamma(1)=x_{i_l}
\quad \text{ for } \quad 
l = 2, \dots, i_{\mathcal{C}_j}
$$ 
its image is a one-chain in 
$Z_1(\{ \tilde{J}_\tau \leq B \}, \{ \tilde{J}_\tau \leq A \})$
with boundary  
$$x_{i_l} - x_{i_1} \in C_0(\{ \tilde{J}_\tau \leq A \}).$$
It turns out that 
\begin{equation}\label{homologies_linearly_independent}
\gamma_{j,2}, \dots, \gamma_{j,\mathcal{C}_j} 
\quad \text{ generate } 
\mathcal{C}_j-1 
\text{  elements of } \quad 
H_1(\{ \tilde{J}_\tau \leq B \}, \{ \tilde{J}_\tau \leq A \}), 
\end{equation}
which are linearly independent. 
To prove \eqref{beta_1_inequality} we show that any
$$\sum_{2 \leq h \leq \mathcal{C}_j} n_h \gamma_{j,h}$$ 
with not all $n_h=0$ cannot be written as 
\begin{equation}\label{eq_homology_contradiction}
\sum_{2 \leq h \leq \mathcal{C}_j} n_h \gamma_{j,h} 
= 
\partial_2 \mathfrak{c}_2 + \mathfrak{c}_1
\end{equation}
with 
$
\mathfrak{c}_2 
\in 
C_2(\{ \tilde{J}_\tau \leq B \}, \{ \tilde{J}_\tau \leq A \})
$
and
$
\mathfrak{c}_1 
\in 
C_1(\{ \tilde{J}_\tau \leq A \}).
$

\

In fact let us apply the boundary operator $\partial_1$ 
to both sides of the latter equation. 
As not all $n_h$ are zero,  
$\partial_1 \left( \sum_h n_h \gamma_{j,h} \right)$ 
is non-trivial in 
$C_0(\{ \tilde{J}_\tau \leq A \})$. 
Clearly 
$\partial_1 \circ \partial_2 = 0$, 
so to achieve \eqref{eq_homology_contradiction} we would need  
$\partial_1 \mathfrak{c}_1$ 
to be in 
$C_0(\{ \tilde{J}_\tau \leq A \})$ 
a non-trivial  linear combination of the points 
$x_{i_1}, \dots, x_{i_l}$. 
However, since 
$x_{i_1}, \dots, x_{i_l}$ 
lie in different components of 
$\{ \tilde{J}_\tau \leq A \}$, 
there is no chain 
$\mathfrak{c}_1 \in C_1(\{ \tilde{J}_\tau \leq A \})$ 
with this property. This shows \eqref{homologies_linearly_independent}.
Repeating this reasoning for every component of $\Xi$  
we obtain
$$
\beta_1 
\geq 
\sum_{j=1}^\mathcal{C} (\mathcal{C}_j - 1) 
= 
p - \mathcal{C},
$$ 
since 
$\sum_{j=1}^\mathcal{C} \mathcal{C}_j = p$. This shows \eqref{beta_1_inequality}.
Now the relative Morse inequalities imply 
$$q = C_1 = C_1 - C_0 \geq \beta_1 - \beta_0 \geq p - \mathcal{C},$$  
contradicting our assumptions. 
\end{proof}
 
 \medskip

In some particular cases we obtain the following  corollary, cf. Theorem \ref{t:pinching} (ii). 

\begin{corollary}
Suppose that $K$ satisfies $\frac{K_{\max}}{K_{\min}} \leq2^{\frac{2}{n-2}}$, 
that it has  $p$ local maxima and $q$ critical points of index $n-1$ with 
negative Laplacian. 
Then \eqref{eq:scin} admits a positive solution provided $q < p-1$. 
\end{corollary}

\begin{proof}
In the  theorem choose the connected set $\Xi = S^n$, see \eqref{eq:unovEm}, 
\eqref{eq:Pm-to-tPm}.  
\end{proof}

\medskip
 
We next state a related result, proved with similar techniques.  
 
\begin{theorem}
Let $(M,g)$ be as in Theorem \ref{t:pq}. 
Suppose $K$ has a local maximum point $z$, and  that there exists a curve $a(t)$ joining $z$ 
to another point $y$ with $$K(y) \geq K(z)$$  such that both the following two properties hold 
\begin{enumerate}[label=(\roman*)]
\item  for all $x_i \neq x_j$ local maxima of $K$ 
$$
\max_t K(a(t))^{\frac{2-n}{n}} 
< 
\left( 
K(x_i)^{\frac{2-n}{2}} + K(x_j)^{\frac{2-n}{2}} 
\right)^{\frac{2}{n}} \};
$$ 
\item critical points of index $n-1$ in the range 
$$[\min_t K(a(t)),K(z)]$$  
have positive Laplacian. 
\end{enumerate}
Then \eqref{eq:scin} has a positive solution. 
\end{theorem}

\begin{proof}
We can construct  a curve $\tilde{a}(t)$ joining $y$ to another maximum point $\tilde{z}$
of $K$ and such that 
$\min_t K(\tilde{a}(t)) = K(y)$. 
Consider then the  composition 
$\hat{a} := a*\tilde{a}$, 
and the test functions 
$\tilde{\varphi}_{x,\tau} $ 
as in \eqref{eq:var-test} for $x$ in the image of the curve $\hat{a}$. 
By Lemma \ref{l:en-est} and construction of $\hat{a}$, 
we have that the image of this curve in $X$ connects two strict local minima 
$u_{\tau,z}$, $u_{\tau,\tilde{z}}$
of $J_\tau$, and the supremum of $J_\tau$ on the image is bounded above by 
$$
\hat{c}_0 (\min_{t \in [0,1]} K(a(t)))^{\frac{2-n}{n}}  
+ 
o_\tau(1) 
+ 
o_{r_0}(1)
.
$$ 

Consider a mountain-pass path between the strict local minima 
$u_{\tau,z}$, $u_{\tau,\tilde{z}}$ 
of $J_\tau$. Assuming that \eqref{eq:scin} has no solutions, 
by the Palais-Smale condition for $J_\tau$ and by the fact that all critical 
points with uniformly bounded energy of $J_\tau$ as described in Theorem \ref{t:ex-multi} 
are non-degenerate, $J_\tau$ must possess a critical point \underline{of index one} 
at a level less or equal to 
$$
\hat{c}_0 (\min_{t \in [0,1]} K(a(t)))^{\frac{2-n}{n}}  
+ 
o_\tau(1) 
+ 
o_{r_0}(1)
.$$ 
Still by Theorem \ref{t:ex-multi} and  condition (i) this critical point must have 
a simple blow-up at  a critical point $p$ of $K$ of index $n-1$ with 
$$K(p) \in [\min_t K(a(t)),K(z)],$$
which is excluded by assumption (ii). 
\end{proof}

\begin{remark}
The latter result improves the pinching condition of 
Theorem 2 in \cite{[BI]} (if compactified from $\R^n$ to $S^n$) for $K$ Morse and satisfying \eqref{eq:nd}, 
namely 
\begin{enumerate}
\item[(j)] 
$K_{\max} < 2^{\frac{2}{n-2}}\min_t K(x(t)) $; 
\item[(jj)]  
critical points in the range 
$[\min_t K(x(t)),K(z))$ 
are local maxima or have positive Laplacian. 
\end{enumerate}
While the strategy in \cite{[BI]} might be possibly used to relax condition (jj), 
an improvement of (j) requires a more careful  analysis of the loss of compactness, 
as done in \cite{MM1} and \cite{M3}. 
\end{remark}

\section{Non-existence results}\label{ss:n-e}

In this section we prove non-existence results on $S^n$ 
for arbitrarily pinched curvature candidates
of prescribed Morse type and with only one critical point with 
negative Laplacian. We  show that the assumptions of Theorem \ref{t:pinching} 
are sharp both in terms of Morse structure and dimension, cf. Remark \ref{r:Struwe34}.

We construct  a sequence of functions $(K_m)_m$ on $S^n$ with only one local maximum, while
all other critical points have positive Laplacian and converge to the south pole. 
We build the $(K_m)_{m}$ in order to preserve a given Morse structure and to maintain uniform 
$C^3$ bounds. 

\noindent
We denote by 
$\mathtt{y}_i$ for $i = 1, \dots, n+1$
the Euclidean coordinate functions on $\R^{n+1}$ restricted to $S^n$
and by $\mathtt{N}, \mathtt{S}$ the north and south poles respectively, i.e. 
$$
\mathtt{N} 
= 
S^n \cap \{ \mathtt{y}_{n+1} = 1 \}
\quad \text{ and } \quad 
\mathtt{S} = S^n \cap \{ \mathtt{y}_{n+1} = - 1 \}.
$$
Finally we let
$$
\pi_{\mathtt{N}} : S^n \setminus \{\mathtt{N}\} \longrightarrow \R^n
\quad \text{ and } \quad
\pi_{\mathtt{S}} : S^n \setminus \{\mathtt{S}\} \longrightarrow \R^n 
$$ 
denote the stereographic projections from the $\mathtt{N},\mathtt{S}$, whose inverse $\pi_{\mathtt{N}}^{-},\pi_{\mathtt{S}}^{-}$ induce coordinate systems 
on $S^{n}\setminus \{\mathtt{N}\},S^{n}\setminus \{\mathtt{S}\}$, to which we will refer as $\pi_{\mathtt{N}}$ and $\pi_{\mathtt{S}}$ coordinates respectively.  

\
 
Recalling our notation in \eqref{eq:MjK} we have the next result, proved in the Appendix. 
  
\begin{proposition}\label{p:Km}
For every Morse function $\tilde{K} : S^n \longrightarrow \R$ with only one local maximum point 
there exists a sequence of positive functions  $(K_m)_m$ such that 
\begin{enumerate}
\item[a)] 
$\mathcal{M}_j(K_m) = \mathcal{M}_j(\tilde{K})$ 
for all $j = 0, \dots, n$ and 
$K_m$ has only one local maximum point at $\mathtt{N}$,  while all other critical points of $K_m$ converge to $\mathtt{S}$; 
\item[b)]  
there exists a neighbourhood 
$U \subseteq S^n$ of $\mathtt{S}$ and $c > 0$ 
such that 
$$\Delta K_m \geq c \quad \text{ on } \quad U;$$  
\item[c)]  
$K_m \longrightarrow K_0$ in $C^3(S^n)$, 
where $K_0$ is a positive monotone non-decreasing function in 
$\mathtt{y}_{n+1}$, 
affine and non-constant in $\mathtt{y}_{n+1}$ outside of a small neighbourhood of $\mathtt{S}$.  
\end{enumerate}
\end{proposition}

\subsection{Uniform bounds away from the poles }
 
We consider the sequence $(K_m)_m$ given by Proposition \ref{p:Km} and
a sequence of positive solutions to 
\begin{equation}\label{eq:sc-m}
L_{g_{S^n}} u_m 
= 
K_m u_m^{\frac{n+2}{n-2}}  
\quad \text{ on } \quad 
(S^n,g_{S^n}). 
\end{equation}
Even without  assuming uniform energy bounds as in Theorem \ref{t:ex-multi},  we aim 
to prove that $(u_m)_m$ stays uniformly bounded on compact sets of 
$S^n \setminus \{\mathtt{N}\}$.

By construction, see the first and final steps in the proof of Proposition \ref{p:Km}, the only 
critical points of $K_0$ are $\mathtt{N}$ and a compact set $K_U \subseteq U$, where 
the Laplacian is positive and bounded away from zero. 
By Corollary 1.4 in \cite{cl1} or Theorem 2 in \cite{cl97}  the sequence $(u_m)_m$ is uniformly bounded in $L^\infty$ on compact sets of  
$$S^n \setminus \{ \mathtt{N} \cup K_U \},$$
since $\vert \n K_m\vert $ is bounded away from zero here, hence we only need to focus  on $K_U$. 
 
For doing this, we cannot directly use known results in the literature due to the degenerating 
behaviour of $(K_m)_m$. However, the proof can be obtained combining the preliminary results 
in Subsection \ref{ss:iiisbu}. It will be harder to understand the blow-up behaviour near the 
north pole $\mathtt{N}$. Before proceeding recall Definition \ref{d:isol-simple}.

 \medskip

\begin{lemma} \label{l:is}
Suppose $(u_m)_m$ solves \eqref{eq:sc-m}.  Then the blow-up points in $U$ are isolated simple. 
\end{lemma}

\begin{proof}
The proof  uses also some argument in Section 8 of \cite{kms}, 
but we have here variable curvature. 
For $0<\d \ll 0$  and $R \gg 1$ let 
$\xi_{1,m}, \dots \xi_{N(u_m),m}$ 
be the points given by Proposition \ref{p:kms}. 
As $(u_m)_m$ is uniformly bounded away from 
$\{\mathtt{N}\} \cup K_U$, 
all $\xi_{i,m}$  will lie in a neighbourhood of 
$\{\mathtt{N}\} \cup K_U$. 
Let us denote by 
$$
\xi_{1,m}, \dots, \xi_{N_m,m} 
\quad \text{ with } \quad 
N_m \leq N(u_m)
$$ 
the points contained in a neighbourhood of $K_U$. 

\
 
\noindent
We may assume that with $c_n$ as in \eqref{eq:conf-change-curv} and in $\pi_{\mathtt{N}}$ coordinates $u_m$ solves 
\begin{equation} \label{eq:eq-um} 
- c_n \Delta u_m 
= 
K_m u_m^{\frac{n+2}{n-2}} 
\quad \text{ in } \quad
B_1(0), 
\end{equation}
where  and we identify $K_m$  with 
$K_{m}\circ \pi_{\mathtt{N}}^{-}$. 
For any $m$ we choose $i \neq j$ such that 
\begin{equation}\label{eq:ij}
\vert \xi_{i,m} - \xi_{j,m}\vert  
= 
\min 
\{ 
\vert \xi_{k,m} - \xi_{l,m}\vert  
\: : \; 
k, l \in \{1, \dots N_m\}, k \neq l 
\}. 
\end{equation}
We let 
$\xi_m = \xi_{i,m}$, $s_m = \frac 12 \vert \xi_{i,m} - \xi_{j,m}\vert$ 
and consider  
\begin{equation}\label{eq:zeta_m}
\zeta_m(x) 
= 
s_m^{\frac{n-2}{2}}u_m \left( s_m x + \xi_{m} \right). 
\end{equation}
By definition of $s_m$ and (iii) in Proposition \ref{p:kms} the sequence $(\zeta_m)_m$ has an isolated  blow-up at zero. 
We will prove next that this blow-up is indeed also isolated simple. 

First,  using the classification result in \cite{cgs}, it is standard to show that there exists 
$R_m \longrightarrow \infty$ 
sufficiently slowly such that 
\begin{equation}\label{eq:conv-um}
	\left\| 
	\zeta_m(0)^{-1} \zeta_m \left( \zeta_m(0)^{-\frac{2}{n-2}} \cdot  \right) 
	- 
	\left( 1 + k_m \vert \cdot \vert ^2 \right)^{\frac{2-n}{2}}  
	\right\| _{C^2(B_{R_m}(\xi_m))} 
	\longrightarrow 
	0,
\end{equation}
where $k_m = \frac{1}{n(n-2)c_n} K_m(\xi_m)$, cf. Proposition 2.1 in \cite{yy}.

Assuming by contradiction that the blow-up of $\zeta_m$ at $0$ is not isolated simple, 
let $\overline{w}_m$ be as in \eqref{eq:barwm} replacing $\overline{u}_m$ by $\overline{\zeta}_m$. 
By \eqref{eq:conv-um} then  $\overline{w}_m$ has a first critical point for $r$ of order 
$\overline{\zeta}_m(0)^{-\frac{2}{n-2}}$ and, if the blow-up of $\zeta_m$ is not isolated simple,  
$$
\tilde{s}_m = \inf  \{ s > R_m \zeta_m(0)^{-\frac{2}{n-2}} \; : \;  \overline{w}'_m(s)=0 \} 
$$  
is well defined and $\tilde{s}_m \ll 1$. If we let $\tilde{\zeta}_m(x) 
= \tilde{s}_m^{\frac{n-2}{2}}\zeta_m \left( \tilde{s}_m x \right)$, 
then  $\tilde{\zeta}_m$ satisfies 
\begin{equation}\label{eq:tildezetam}
	- c_n \Delta \tilde{\zeta}_m = \tilde{K}_m(x) \tilde{\zeta}_m^{\frac{n+2}{n-2}}; 
	\qquad \tilde{K}_m (x) = K_m(\xi_m + \hat{s}_m x), \quad \hat{s}_m = s_m \tilde{s}_m,  
\end{equation}
and has an isolated blow-up at zero. From Lemma \ref{l:sing-prof} we deduce  
$$
\tilde{\zeta}_m(0) \tilde{\zeta}_m(x)
\longrightarrow 
a \vert x\vert ^{2-n} + h(x) \geq 0
\quad \text{ in } \quad C^2_{\text{loc}} (\R^n \setminus \{0\}),$$ 
where $h$ is harmonic on $\R^n$ and $a > 0$. 
By the first observation after Lemma \ref{l:sing-prof} the function  $h$ must be constant, 
and passing to the limit for the condition 
$\overline{w}'_m(\tilde{s}_m) = 0$ 
one  finds that  $h \equiv a > 0$, as for (3.4) in \cite{yy}.

%
%
%
%
%

From Lemma \ref{l:sing-prof} and, since $\tilde{\zeta}_m$ has an isolated simple blow-up, it follows that 
\begin{equation}\label{eq:ineq-um-s-p}
\tilde{\zeta}_m(x) 
\leq 
C \tilde{\zeta}_m(0)^{-1} \vert x \vert ^{2-n} 
\quad \text{ for } \quad 
\vert x \vert  
\in 
[R_m \tilde{\zeta}_m(0)^{-\frac{2}{n-2}}, 1]. 
\end{equation}
For $\delta > 0$ fixed, we now let $B_{\delta} := B_{\delta}(0)$, and 
for all $i= 1, \dots, n$ we clearly have  
\begin{equation}\label{eq:eq-1}
\begin{split}
\frac{1}{2^*} \int_{B_{\delta}} \frac{\pa \tilde{K}_m}{\pa x_i} \tilde{\zeta}_m^{2^*} dx 
= &
\frac{1}{2^*} \int_{B_{\delta}} \frac{\pa \tilde{K}_m}{\pa x_i}(0) \tilde{\zeta}_m^{2^*} dx \\
& + 
\frac{1}{2^*} \int_{B_{\delta}} \left( \frac{\pa \tilde{K}_m}{\pa x_i} 
-
\frac{\pa \tilde{K}_m}{\pa x_i}(0) \right) \tilde{\zeta}_m^{2^*} dx. 
\end{split}
\end{equation}
By the uniform $C^3$-bounds on $(K_m)$, see Proposition \ref{p:Km}, the convergence in \eqref{eq:conv-um}, the 
upper bound in \eqref{eq:ineq-um-s-p}, a cancellation by oddness and a change of variables we find that the 
last term in \eqref{eq:eq-1} is of order $o(\hat{s}_m \tilde{\zeta}_m(0)^{-\frac{2}{n-2}})$, so
$$
\frac{1}{2^*} \int_{B_{\delta}} \frac{\pa \tilde{K}_m}{\pa x_i}(0) \tilde{\zeta}_m^{2^*} dx 
= 
\frac{1}{2^*} \int_{B_{\delta}} \frac{\pa \tilde{K}_m}{\pa x_i} \tilde{\zeta}_m^{2^*} dx
+ 
o(\hat{s}_m \tilde{\zeta}_m(0)^{-\frac{2}{n-2}}). 
$$ 
By elliptic regularity theory the upper bound \eqref{eq:ineq-um-s-p} implies  
$$
\vert \n \tilde{\zeta}_m(x)\vert  
\leq 
C \tilde{\zeta}_m(0)^{-1} 
\quad \text{ on } \quad 
\pa B_{\delta}.
$$
 Therefore, from \eqref{eq:poho-trans}
we deduce
$$
\frac{1}{2^*} \int_{B_{\delta}} \frac{\pa \tilde{K}_m}{\pa x_i} \tilde{\zeta}_m^{2^*} dx  
= 
\oint_{\pa B_{\delta}} O( \tilde{\zeta}_m^{2^*} + \vert \n \tilde{\zeta}_m\vert ^2) d \s  
= 
O_\delta (\tilde{\zeta}_m(0)^{-2} ). 
$$
It follows from the last two formulas  that 
\begin{equation}\label{eq:estgradKai}
\frac{\pa \tilde{K}_m}{\pa x_i}(0)  
=  
O_\delta ( \tilde{\zeta}_m(0)^{-2}  ) 
+ 
o(\hat{s}_m \tilde{\zeta}_m(0)^{-\frac{2}{n-2}}). 
\end{equation} 
We next rewrite \eqref{eq:poho} for $\tilde{\zeta}_m$ as 
\begin{equation}\label{eq:poho-2222} 
\begin{split}
\frac{1}{2^*} \int_{B_{\delta}} \sum_i x_i  & \frac{\pa \tilde{K}_m}{\pa x_i}(0) \tilde{\zeta}_m^{2^*} dx   +   
\frac{1}{2^*}  \int_{B_{\delta}} \sum_i x_i 
\left(  
\frac{\pa \tilde{K}_m}{\pa x_i} 
- 
\frac{\pa \tilde{K}_m}{\pa x_i}(0) 
\right)  
\tilde{\zeta}_m^{2^*} dx 
\\ 
& - 
\frac{1}{2 \cdot 2^*} 
\oint_{\pa B_{\delta}} \tilde{K}_m \tilde{\zeta}_m^{2^*} d \s 
=  
c_n \oint_{\pa B_{\delta}} B(1/2, x, \tilde{\zeta}_m, \n \tilde{\zeta}_m) d \s.  
\end{split}	
\end{equation}
Using the same reasoning as after \eqref{eq:eq-1}, one finds that 
$$
\int_{B_{\delta}}  x_i  \tilde{\zeta}_m^{2^*} dx  
= 
o ( \tilde{\zeta}_m(0)^{-\frac{2}{n-2}} ).
$$
From these formulas and \eqref{eq:estgradKai} we then deduce that 
\begin{equation}\label{eq:ee2}
\int_{B_{\delta}} \sum_i x_i \frac{\pa \tilde{K}_m}{\pa x_i}(0) \tilde{\zeta}_m^{2^*} dx
= 
o_\delta ( \tilde{\zeta}_m(0)^{-\frac{2(n-1)}{n-2}} ) 
+ 
o ( \hat{s}_m \tilde{\zeta}_m(0)^{-\frac{4}{n-2}} ). 
\end{equation}
Still using the uniform $C^3$-bounds on $(K_m)$, the convergence in \eqref{eq:conv-um}, the 
upper bound in \eqref{eq:eq-1} and a change of variables we find that  
with some $l_n>0$
\begin{equation}\label{eq:ee3}
\int_{B_{\delta}} \sum_i x_i 
\left(  
\frac{\pa \tilde{K}_m}{\pa x_i} 
- 
\frac{\pa \tilde{K}_m}{\pa x_i}(0) 
\right)  
\tilde{\zeta}_m^{2^*} dx 
= 
l_n \hat{s}_m ( \Delta K_m (\xi_m) + o_m(1) ) \tilde{\zeta}_m(0)^{-\frac{4}{n-2}}. 
\end{equation}
Moreover, since 
$
\tilde{\zeta}_m(x) 
\leq 
C \tilde{\zeta}_m(0)^{-1} \vert x \vert ^{2-n} 
$ 
on $\pa B_{\delta}$, we have 
$$
\frac{1}{2 \cdot 2^*} \oint_{\pa B_{\delta}} \tilde{K}_m \tilde{\zeta}_m^{2^*} d \s 
= 
O_\delta ( \tilde{\zeta}_m(0)^{-2^*} ), 
$$
so recalling \eqref{eq:lim-B} we get from \eqref{eq:poho-2222} and the latter estimates that, for $\delta$ small 
\begin{equation*}
\begin{split}
\frac{l_n}{2^*} \hat{s}_m ( \Delta K_m (\xi_m) & + o_m(1) ) \tilde{\zeta}_m(0)^{-\frac{4}{n-2}}  \\
& +  
\frac{(n-2)^2}{2} h(0) \omega_n \frac{c_n + o_m(1)}{ \tilde{\zeta}_m(0)^2} 
= 
o_\delta ( \tilde{\zeta}_m(0)^{-\frac{2(n-1)}{n-2}} ), 
\end{split}
\end{equation*}
a contradiction to $h(0) = a > 0$  and the fact that 
$\Delta K_m (\xi_m) $ 
is positively bounded away from zero. We  hence proved that $\zeta_m$ 
has an isolated simple blow-up at zero.

 \medskip
 
The  exactly  same strategy, but using the second observation 
after Lemma \ref{l:sing-prof},  then shows
 \begin{equation}\label{claim}
 	2 s_m = \vert \xi_{i,m} - \xi_{j,m}\vert  
 	\not\to 
 	0
 	\quad \text{ as }\quad
 	m \longrightarrow \infty, 
 \end{equation} 
as for Section 8 in \cite{kms},
 proving that the blow-ups of $u_m$ in $U$ are isolated. 
 Repeating once more the argument used above for $\tilde{\zeta}_m$ 
   shows that the blow-ups of $u_m$ in $U$ are indeed 
 also isolated simple, which is the desired result. 
\end{proof}
 
 \medskip 
 
\begin{proposition}\label{p:u-b}
For $(K_m)_m$ given by Proposition \ref{p:Km} let  $u_m$ solve \eqref{eq:sc-m} with 
$n \geq 5$. 
Then $(u_m)_m$ is uniformly bounded on the compact sets of $S^n \setminus \{\mathtt{N}\}$. 
\end{proposition}

\begin{proof}
Using the notation in the previous proof, it is sufficient to prove that no blow-up occurs at points in $K_U$. 
We know by Lemma \ref{l:is} that such blow-ups would be isolated simple 
and therefore they could be at most finitely-many. Let $\xi_m \longrightarrow \xi_U$ be a blow-up point in $K_U$. 
Then by Lemma \ref{l:sing-prof} and the Harnack inequality we find that
in $\pi_{\mathtt{N}}$ coordinates 
\begin{equation*}
u_m(\xi_m) u_m(x+\xi_m) \longrightarrow a \vert x\vert ^{2-n} + h(x)
\quad \text{ in } \quad 
C^2_{\text{loc}}(\R^n \setminus S),
\end{equation*}
where $S$ is a finite set, $a > 0$ and $h$  harmonic near $0\in S$.
Moreover $h(0) \geq 0$, see the comments after Lemma \ref{l:sing-prof}.
By Lemma \ref{l:sing-prof} there exists some \underline{fixed} $r > 0$ so that 
the upper bound \eqref{eq:ub-um} holds on $\pa B_{r/2}(0)$. 
Hence and by \eqref{eq:lim-B} we obtain
$$
\frac{r}{2} \oint_{\pa B_{r/2}(\xi_m)} K_m u_m^{2^*} d\s 
=
\frac{O(1)}{u_m(\xi_m)^{2^*}}
$$
and
$$
\oint_{\pa B_{r/2}(\xi_m)} B(r/2, x, u_m, \n u_m) d \s 
\leq
\frac{o_m(1)}{u_m(\xi_m)^{2}}.
$$
Moreover, reasoning as for \eqref{eq:ee2} and \eqref{eq:ee3}, but on a ball of fixed radius, we find that for some $l_n > 0$ 
$$
\int_{B_{r/2}(\xi_m)} \sum_i x_i \frac{\pa K_m}{\pa x_i} u_m^{2^*} dx 
= 
\frac
{l_n \D K_m(\xi_m) + o_m(1)}
{u_m(\xi_m)^{\frac{4}{n-2}}}, 
$$
which immediately leads to a contradiction to \eqref{eq:poho}, since $n\geq 5$ and 
$$\D K_m(\xi_m)  \geq c/2 > 0.$$  
This concludes the proof. 
\end{proof}

 \subsection{Conclusion}

Here we prove our non-existence result, Theorem \ref{t:non-ex}, showing that sequences of solutions to \eqref{eq:sc-m} 
can neither have a non-zero limit nor develop blow-ups, which is impossible.

\begin{lemma}\label{l:no-sing}
Let $K_0$ be a monotone function as in Proposition \ref{p:Km}. Then neither 
\begin{equation}\label{eq:regreg}
L_{g_{S^n}} u = K_0 u^{\frac{n+2}{n-2}} 
\quad \text{ on } \quad
S^n, 
\end{equation}
nor 
\begin{equation}\label{eq:singsing}
L_{g_{S^n}} u = K_0 u^{\frac{n+2}{n-2}} 
\quad \text{ on } \quad
S^n \setminus \{ \mathtt{N} \} 
\end{equation}
admits positive solutions.
\end{lemma}

\begin{proof}
Non existence for \eqref{eq:regreg} simply follows from the Kazdan-Warner obstruction. 
Arguing by contradiction for \eqref{eq:singsing}, 
we obtain in  $\pi_{\mathtt{S}}$ coordinates and by conformal invariance of the equation a positive solution $u$ of the problem  
\begin{equation} \label{eq:change}
- c_n \Delta u 
= 
K_0 u^{\frac{n+2}{n-2}}  
\quad \text{ in } \quad
\R^n \setminus \{0\},  
\end{equation}
where we are identifying $K_0$ with $K_0 \circ \pi_{\mathtt{S}}^{-1}$, which is radially non-increasing and somewhere strictly decreasing.
Since the solution of \eqref{eq:singsing}
is smooth near $\mathtt{S}$, the solution $u$ of \eqref{eq:change} satisfies 
\begin{equation}\label{eq:decdec}
u(x) \leq C \vert x\vert ^{2-n}
\quad \text{ and } \quad 
\vert \nabla u(x)\vert  \leq C \vert x\vert ^{1-n} 
\quad \text{ for } \quad
\vert x\vert  \longrightarrow \infty 
\end{equation}
for some positive and fixed constant $C$. Let us write  the Pohozaev identity 
in the complement of a ball, i.e. on 
$$A_\e := \R^n \setminus B_\e(0).$$ 
By \eqref{eq:decdec} no boundary terms at infinity are involved, whence
\begin{equation}\label{eq:poho2}
\frac{1}{2^*}\int_{A_\e} u^{2^*} \sum_i x_i \frac{\pa K_0}{\pa x_i}  \,  dx 
= 
\frac{1}{2^*} \oint_{\pa A_\e} \langle x, \nu \rangle K_0  u^{2^*} d \s 
+ 
c_n \oint_{\pa A_\e} B(\e, x, u, \n u) \, d \s, 
\end{equation}
see \eqref{eq:poho} and the subsequent formula. By Theorem 1.1 in \cite{zha02}  
\begin{equation}\label{eq:ubuuu} 
\; \exists \; C>0 \; : \; u(x) \leq C \vert x\vert ^{\frac{2-n}{2}}  
\quad \text{ as } \quad 
0 \neq x \longrightarrow 0. 
\end{equation} 
We now consider two cases. 

\medskip 

\noindent {\bf Case 1.} There exists $C > 0$ such that 
$$
C^{-1} \vert x\vert ^{\frac{2-n}{2}} 
\leq 
u(x) 
\quad \text{ as } \quad 
0 \neq x \longrightarrow 0. 
$$
In this case there exists by Theorem 1 in \cite{tz} a singular, radial Fowler's solution 
$$
- \Delta u_0 
= 
\kappa \, u_0^{\frac{n+2}{n-2}} 
\quad \text{ in } \quad
\R^n \setminus \{0\}
\quad \text{ with } \quad 
\kappa = c_n^{-1} K_0(\mathtt{N})
$$
with negative Hamiltonian energy, cf. Subsection \ref{ss:sing-sol}, 
such that 
$$
u(x) = (1 + O(\vert x\vert ^2)) u_0(x).
$$  
Since the unit normal to $A_\e$ points toward the origin, 
the right-hand side of \eqref{eq:poho2} 
is by Lemma \ref{l:cons-law} positive for $\e$ sufficiently small. 
On the other hand the left-hand side of \eqref{eq:poho2} 
is negative by radial monotonicity of $K_0 \circ \pi^{-1}$ and  positivity of $u$, 
so we reach a contradiction.

\medskip 
 
\noindent {\bf Case 2.} Suppose there exists $x_m \longrightarrow 0$ such that 
\begin{equation}\label{eq:um00}
u(x_m) 
= 
o_m(1) \vert x_m\vert ^{\frac{2-n}{2}}. 
\end{equation}
The upper bound in \eqref{eq:ubuuu} yields a Harnack inequality 
for $u$ on annuli of the type $B_{2s}(0) \setminus B_{s/2}(0)$, cf. the 
proof of Lemma 2.1 in \cite{yy}. Thus  by elliptic regularity theory 
there exists $\e_m \searrow 0$ such that for
$x \in B_{2 \e_m}(0) \setminus B_{\e_m/2}(0)$ 
$$
u(x) 
=
o_m(1) \vert \e_m\vert ^{\frac{2-n}{2}}
\quad \text{ and } \quad 
\vert \n u(x)\vert  = o_m(1) \vert \e_m\vert ^{-\frac{n}{2}} 
. 
$$
This and \eqref{eq:poho} imply that for such an $(\e_m)_m$
$$
\frac{1}{2^*} \oint_{\pa A_{\e_m}} \langle x, \nu \rangle K_0(x)  u^{2^*} d \s 
+ 
c_n \oint_{\pa A_{\e_m}} B(\e_m, x, u, \n u) \, d \s 
\longrightarrow 
0, 
$$
contradicting \eqref{eq:poho2} as in the previous case. 
\end{proof}
 
 \medskip 
 
\noindent As an immediate consequence of Proposition \ref{p:u-b} and Lemma \ref{l:no-sing} we obtain the following result. 
 
\begin{proposition}\label{p:zwl}
For $(K_m)_m$ as in Proposition \ref{p:Km} let $u_m > 0$ solve \eqref{eq:sc-m} with $n \geq 5$.
Then $(u_m)_m$ converges to zero in $C^2_{\text{loc}}(S^n \setminus \{\mathtt{N}\})$. 
\end{proposition}

\medskip \noindent
We next analyse also the case of zero-limit  in $C^2_{\text{loc}}(S^n \setminus \{\mathtt{N}\})$, showing that a non-zero one can be obtained after a proper dilation.

\begin{lemma}\label{l:non-trivial}
Let $(u_m)_m$ be as in Proposition \ref{p:zwl}.  
Then, writing \eqref{eq:sc-m} in $\pi_{\mathtt{S}}$ coordinates, i.e.
\begin{equation}\label{eq:Km}
- \D  u_m 
= 
K_m u_m^{\frac{n+2}{n-2}} 
\quad \text{ in } \quad
\R^n, 
\end{equation}
there is near the north pole $\mathtt{N}$ a blow-down $(v_m)_m$ of $(u_m)_m$ of the form 
\begin{equation}\label{eq:scale-vm}
v_m(x) 
= 
\mu_m^{\frac{n-2}{2}} u_m(\mu_m x)
\quad \text{ with } \quad 
\mu_m 
\longrightarrow 
0, 
\end{equation}
such that up to a subsequence $(v_m)_m$ has a non zero limit in $C^2_{\text{loc}}(\R^n \setminus \{0\})$.  
\end{lemma}
 
\begin{proof}
We blow-up the metric $g_{S^n}$ conformally near $\mathtt{N}$ in order to obtain a metric 
$$
\tilde{g} = \tilde{u}^{\frac{4}{n-2}} g_{S^n}
\quad \text{ with } \quad 
\tilde{u} \simeq \vert x\vert ^{\frac{2-n}{2}}
\quad \text{ near } \quad 
x = 0$$ 
in the above coordinates and  
with a cylindrical end and  bounded geometry. If 
$$\tilde{u}_m = \tilde{u}^{-1} u_m,$$ then by \eqref{eq:covar} $\tilde{u}_m$ satisfies 
$$
L_{\tilde{g}} \tilde{u}_m 
= 
K_m \tilde{u}_m^{\frac{n+2}{n-2}} 
\quad \text{ on } \quad
(S^n \setminus \{ \mathtt{N} \},\tilde{g}). 
$$
By (1.7) in \cite{cl2} we have 
$
u_m(x) 
\leq 
C \vert x\vert ^{\frac{2-n}{2}}
$, 
whence $(\tilde{u}_m)_{m}$ is uniformly bounded.
Note that the dilation  in  \eqref{eq:scale-vm} corresponds to a translation along the cylindrical end in the metric $\tilde{g}$ and yields  
$v_m(x) \leq C \vert x\vert ^{\frac{2-n}{2}}$. 
 
\
Using the assumption on the zero-limit in $C^2_{\text{loc}}$ of $u_m$ on 
$S^{n} \setminus \{\mathtt{N}\}$, 
elliptic regularity theory and the uniform bound on $\tilde{u}_m$,
and arguing by contradiction 
$$
v_m \longrightarrow 0
\quad \text{ in } \quad C^2_{\text{loc}}(\R^n \setminus \{0\})
\quad \text{ for every choice of } \quad 
\mu_m \searrow 0$$ 
would imply 
$\tilde{u}_m \longrightarrow 0$  uniformly on $S^n \setminus \{\mathtt{N}\}$. 
We then use elliptic estimates for
$$
- 4 \frac{n-1}{n-2} \Delta_{\tilde{g}} \tilde{u}_m 
+ 
R_{\tilde{g}} \tilde{u}_m  
= 
K_m \tilde{u}_m^{\frac{n+2}{n-2}}
\quad \text{ on } \quad
(S^n \setminus \{ \mathtt{N} \},\tilde{g})
$$
to show, that for $x$ in the cylindrical end of 
$(S^{n} \setminus \{\mathtt{N}\}, \tilde{g})$, 
 where $R_{\tilde{g}}$ is positive,  
$$
\| \tilde{u}_m \|_{L^\infty(B_1(x))} 
\leq 
C \| \tilde{u}_m \|_{L^\infty(B_1(x))}^{\frac{n+2}{n-2}}. 
$$
Here the metric ball around $x$ is taken with respect to  $\tilde{g}$. 
Since the latter norm tends to zero for $m \longrightarrow \infty$, $\tilde{u}_m$
must be identically zero for $m$ large near the cylindrical end, 
contradicting the  positivity  of $u_m$. 
\end{proof}
 
 \medskip

We next perform a blow-down as in Lemma \ref{l:non-trivial} at 
{\em slowest possible rate}, i.e. working in $\pi_{\mathtt{S}}$ coordinates
we can choose, e.g. with a concentration-compactness argument,  $\bar \mu_m \searrow 0$ with the properties 
\begin{enumerate}
\item \quad
$
\bar{v}_m(x) 
= 
\bar \mu_m^{\frac{n-2}{2}} u_m(\bar \mu_m x)
$ 
converges in 
$C^2_{\text{loc}}(\R^n \setminus \{0\})$ 
to a non-zero limit;
\item \quad 
if $\frac{\hat{\mu}_m}{\bar \mu_m} \longrightarrow 0$, then 
$\hat \mu_m^{\frac{n-2}{2}} u_m(\hat \mu_m x)$ 
converges to zero in $C^2_{\text{loc}}(\R^n \setminus \{0\})$. 
\end{enumerate}
 
\begin{lemma}\label{l:reg-bubble-limit}
Up to a subsequence $(\bar v_{m})_{m}$ converges in $C^2_{\text{loc}}(\R^n \setminus \{0\})$   to a regular bubble.
\end{lemma}

\begin{proof}
If $v_0$ is the limit of 
$\bar{v}_m$ in $C^2_{\text{loc}}(\R^n \setminus \{0\})$, 
it satisfies 
$$
- \Delta v_0 
= 
\kappa \, v_0^{\frac{n+2}{n-2}} 
\quad \text{ in } \quad
\R^n \setminus \{0\}, 
\quad \text{ where } \quad 
\kappa = c_n^{-1} K_0(\mathtt{N}). 
$$
Due to the classification result in Corollary 8.2 of \cite{cgs} 
we need to prove that $v_0$ has a removable singularity near zero. 
Assume by contradiction that $v_0$ is singular there. Then $v_0$ must be radially symmetric
by Theorem 8.1 in \cite{cgs}. 
Singular radial solutions are  classified as described in Subsection \ref{ss:sing-sol} as Fowler's solutions and by positivity of
$v_H$ for any such solution there exists $c > 0$ such that
$$v_0 \geq \frac{c}{\vert x\vert ^{\frac{n-2}{2}}}.$$ 
Hence we proved that in case of a singular limit $v_0$, 
$$
\bar{v}_m \longrightarrow v_0 
\quad \text{ in } \quad
C^2_{\text{loc}}(\R^n \setminus \{0\})
\quad \text{ and } \quad 
v_0(x) \geq \frac{c}{\vert x\vert ^{\frac{n-2}{2}}}, 
$$
which would violate the above condition (ii) on $\bar{\mu}_m$. This concludes the proof. 
\end{proof}

 \medskip
 
\begin{lemma}\label{l:outside}
If $(\bar v_m)_m$ is  as above, then there exists $C > 0$ such that 
$$
u_m(x) 
\leq C 
\bar{\mu}_m^{\frac{n-2}{2}} d_{S^n}(x,\mathtt{N})^{2-n}
\quad \text{ for } \quad 
d_{S^n}(x,\mathtt{N}) \geq \bar \mu_m. 
$$
\end{lemma}

 \medskip
The lemma is proved in the appendix. We next consider a Kelvin inversion around a sphere of radius 
$\check{\mu}_m \longrightarrow 0$ 
with 
$\frac{\check{\mu}_m}{\bar \mu_m} \longrightarrow 0$. 
In $\pi_{\mathtt{S}}$ stereographic coordinates this corresponds to the map 
$$
x \mapsto \frac{\check{\mu}_m^2 x }{\vert x\vert ^2}. 
$$
Letting
\begin{equation}\label{eq:checku_m}
\check{u}_m (x) 
= 
\frac{\check{\mu}_m^{n-2}}{\vert x\vert ^{n-2}} 
u_m \left( \frac{\check{\mu}_m^2 x }{\vert x\vert ^2} \right), 
\end{equation}
we obtain from \eqref{eq:Km} a sequence of functions $\check{u}_m$ satisfying 
\begin{equation}\label{eq:sat}
-  c_n \Delta \check{u}_m
= 
\check{K}_m \check{u}_m^{\frac{n+2}{n-2}}
\quad \text{ in } \quad
B_1(0),
\quad \text{ where } \quad 
\check{K}_m(x)  = K_m \left( \frac{\check{\mu}_m^2 x }{\vert x\vert ^2} \right).
\end{equation}
As the functions $\check{K}_m$ are highly oscillating near $x=0$, we 
lose uniform Lipschitz bounds compared to $(K_m)_{m}$. 
More precisely, let $\mathring{K}_m$ denote the functions $K_m$ 
reflected with respect to the hyperplane $\{\mathtt{y_{n+1}} = 0\}$ in $\R^{n+1}$. 
By direct calculation 
$\check{K}_m(x) = \mathring{K}_m (\check{\mu}_m^{-2} x)$ 
for $x \in B_1(0)$, where we are indentifying $\mathring{K}_m$ 
with $\mathring{K}_m\circ \pi^{-}_{\mathtt{S}}$ as before.
This implies 
\begin{equation}\label{eq:bd-grad-check-K-m}
\vert \n \check{K}_m(x)\vert  
\leq 
\frac{C}{\check{\mu}_m^2}
\quad \text{ for } \quad 
x \in 
B_1(0). 
\end{equation}
However, since 
$$
K_0(x) 
= 
\kappa 
- 
(\kappa_0 + o_m(1) ) \vert x\vert ^2 
+ 
O(\vert x\vert ^3)
\quad \text{ for } \quad 
\vert x\vert  \leq \d
$$
and some 
$\kappa_0 > 0$ 
by (c) of Proposition \ref{p:Km}, we have
\begin{equation}\label{eq:check-K-m}
\check{K}_m(x) 
= 
\kappa 
- 
\kappa_0 (1 + o_m(1)) \frac{\check{\mu}_m^4}{\vert x\vert ^2} 
+ 
O(\check{\mu}_m^6 \vert x\vert ^{-3})
\quad \text{ for } \quad 
\vert x\vert  \geq \frac{\check{\mu}_m^2}{\d}. 
\end{equation}
Let $U_0$ be as in \eqref{eq:U0} and define 
$$U_{a,\l}(x) = \lambda^{\frac{n-2}{2}} U_0(\lambda(x-a))$$ 
for $a \in \R^n$ and $\l > 0$. 
By Lemma \ref{l:reg-bubble-limit} then $u_m$ is on a proper annulus centred at 
$x = 0$  close in $W^{1,2}$ to a multiple, which depends on $K_0(\mathtt{N})$, of 
$U_{a_m,\lambda_m}$ with  $\lambda_m \simeq \bar{\mu}_m^{-1}$. 
As
$$u_m(x) \leq C \vert x\vert ^{\frac{2-n}{2}}$$ by (1.7) in \cite{cl2}, we find that $\lambda_m \vert a_m\vert $ is uniformly bounded. 
By direct computation the inversion in \eqref{eq:checku_m} sends 
$U_{a_m,\l_m}$ 
into 
$U_{\check{a}_m,\check{\l}_m}$, 
where 
$$
\check{a}_m 
= \lambda_m^2 \check{\mu}_m^2 \frac{a_m}{1+ \lambda_m^2 \vert a_m\vert ^2}
\quad \text{ and } \quad 
\check{\l}_m = \frac{1+ \lambda_m^2 \vert a_m\vert ^2}{\lambda_m \check{\mu}_m^2}. 
$$
Note, that $\check{\l}_m \vert \check{a}_m\vert $ is uniformly bounded, as  $\lambda_m \vert a_m\vert $ is. 
Hence  
$$
\; \exists\; y_m \longrightarrow 0 \; : \; 
\check{u}_m(y_m) 
\simeq 
\left( \frac{\check{\mu}_m^2}{\bar{\mu}_m} \right)^{\frac{2-n}{2}} 
\longrightarrow 
\infty
$$ 
and $\check{u}_m$ develops a bubble at a scale
$$\frac{\check{\mu}_m^2}{\bar{\mu}_m} \longrightarrow 0.$$
Since the Kelvin inversion and the above bound on $u_m$ yield  the condition  
$$\check{u}_m(x) \leq C \vert x\vert ^{\frac{2-n}{2}},$$ 
$x = 0$ is the only blow-up point for $(\check{u}_m)_m$. 
Moreover by Lemma \ref{l:outside} we also deduce
$$
\max \check{u}_m 
\simeq 
\left( \frac{\check{\mu}_m^2}{\bar{\mu}_m} \right)^{\frac{2-n}{2}}
.
$$ 
Note that from the regular bubbling profile, cf. Lemma \ref{l:reg-bubble-limit}, the radial average 
$$
\bar{\check{w}}_m(r)
= 
r^{\frac{n-2}{2}} \intbar_{\pa B_r(x_m)} \check{u}_m d \s
$$
has a unique critical point for $r$ of order $\frac{\check{\mu}_m^2}{\bar{\mu}_m}$, 
see \eqref{eq:barwm}.  
If there is another critical point at some $\check{r}_m \longrightarrow 0$, it must be $\check{r}_m \gg \frac{\check{\mu}_m^2}{\bar{\mu}_m}$. 
Therefore we can choose $\check{\mu}_m$ so that 
$\bar{\check{w}}_m$ has a unique critical point in 
$
\left[ 
\frac{\check{\mu}_m^2}{\bar{\mu}_m} 
, 
1 
\right]
$. 
Despite the oscillations of the $\check{K}_m$'s we have  the following result, also proven in the appendix. 

\begin{lemma}\label{l:sing-prof-anyway}
Suppose that  
$\check{\mu}_m \ll \bar{\mu}_m$ 
is chosen so that 
$\bar{\check{w}}_m$ 
has a unique critical point in 
$
\left[ 
\frac{\check{\mu}_m^2}{\bar{\mu}_m} 
, 
1 
\right]
$. 
Then the same conclusions of Lemma \ref{l:sing-prof} hold true. 
\end{lemma}

\medskip
 
\noindent  
We can finally prove our non-existence result, yielding also Theorem \ref{t:non-ex}. 
 
\begin{theorem}
Suppose that $(K_m)_m$ is as in Proposition \ref{p:Km}. Then for $m$ large 
problem \eqref{eq:sc-m} has no positive solutions. 
\end{theorem}

\begin{proof}
Assume by contradiction that \eqref{eq:sc-m} possesses positive solutions for all $m$. 
We saw in Proposition \ref{p:u-b} that $(u_m)_m$ is uniformly bounded on 
$S^n \setminus \{\mathtt{N}\}$, 
so up to a subsequence we have that 
$$
u_m \longrightarrow u_0 
\quad \text{ in } \quad 
C^2_{\text{loc}}(S^n \setminus \{\mathtt{N}\}), 
$$
where $u_0$ solves 
$$
L_{g_{S^n}} u_0  
= 
K_0 \, u_0^{\frac{n+2}{n-2}} 
\quad \text{ on } \quad 
S^n \setminus \{ \mathtt{N}\}
\quad \text{ with } \quad 
K_0 = \lim_m K_m. 
$$
By Lemma \ref{l:no-sing}, $u_0$ can be  neither a regular nor a positive singular solution. 
Therefore we must have 
$u_0 \equiv 0$ 
and can hence apply Lemmas \ref{l:non-trivial} 
and \ref{l:reg-bubble-limit}, letting $\bar{\mu}_m$ as in Lemma \ref{l:reg-bubble-limit}. 
 
Working in $\pi_{\mathtt{S}}$ coordinates and choosing $\check{\mu}_m$ properly, 
$\check{u}_m$ defined in \eqref{eq:checku_m} satisfies the assumptions of Lemma \ref{l:sing-prof-anyway}.
Therefore we have for $(\check{u}_m)_m$ the conclusion of Lemma \ref{l:sing-prof}. 
Let as before $y_m$ be a global maximum of $\check{u}_m$. 
As remarked after Lemma \ref{l:sing-prof}, we have that 
$$
\check{u}_m(y_m) \check{u}_m 
\longrightarrow 
a \vert x\vert ^{2-n} + h(y)
\quad \text{ in } \quad 
C^2_{\text{loc}}(\R^n \setminus \{0\}),
$$
where $a > 0$ and $h \geq 0$ is identically constant. 
From this and  \eqref{eq:lim-B} we find
\begin{equation}\label{eq:poho-partial}
\oint_{\pa B_1} \check{K}_m \check{u}_m^{2^*} d \s 
= 
o \left( \frac{\check \mu_m^2}{\bar \mu_m} \right)^2
\; \text{ and } \;
\oint_{\pa B_1} B(\rho, x, \check{u}_m, \n \check{u}_m) d \s 
=
o \left( \frac{\check \mu_m^2}{\bar \mu_m} \right)^2. 
\end{equation}
Letting now $\d$ as in \eqref{eq:check-K-m}, from Lemma \ref{l:outside} we find
$$
\check{u}_m 
\leq 
C \left( \frac{\check{\mu}_m^2}{\bar{\mu}_m} \right)^{\frac{2-n}{2}}
\quad \text{ for } \quad 
\vert x\vert  \leq \frac{\check{\mu}_m^2}{\d}.
$$ 
Hence by \eqref{eq:bd-grad-check-K-m}, \eqref{eq:check-K-m} and, as
$\check{\mu}_m$ develops a bubble at scale 
$\frac{\check{\mu}_m^2}{\bar{\mu}_m}$, 
$$
\int_{B_{\frac{\check{\mu}_m^2}{\d}}} 
\sum_i x_i \frac{\pa \check{K}_m}{\pa x_i} \check{u}_m^{2^*} 
dx 
= 
O(\bar{\mu}_m^n)
\quad \text{ and } \quad 
\int_{B_1 \setminus B_{\frac{\check{\mu}_m^2}{\d}}} \sum_i x_i \frac{\pa \check{K}_m}{\pa x_i} \check{u}_m^{2^*} dx 
\geq  
c \bar \mu_m^2, 
$$
where $c > 0$, cf. the discussion after \eqref{eq:check-K-m}. From this we deduce
$$
\int_{B_1}\sum_i x_i \frac{\pa \check{K}_m}{\pa x_i} \check{u}_m^{2^*} dx 
\geq  
c\bar \mu_m^2
,$$ 
yielding a contradiction together with 
\eqref{eq:poho}, \eqref{eq:poho-partial} and
$\frac{\check{\mu}_m}{\bar \mu_m} \longrightarrow 0$. 
%
%
%
 \end{proof}

\begin{remark}\label{r:Struwe34}
In \cite{str05} a non-existence result was proved on $S^2$ for curvature functions 
that are not monotone with respect to any Euclidean coordinate in $\R^3$ restricted to the unit sphere. 
Such functions have two maxima and one saddle point close to  the north pole and in addition 
one non-degenerate minimum near the south pole, hence they are 
{\em reversed} compared to the ones considered in this section. 
 
The proof of the above result in \cite{str05} relies on showing that solutions would be close to a single bubble: in this 
way the left-hand side in \eqref{eq:kw} can be made {\em quantitatively} non-zero (depending on the 
concentration rate of the bubble), even if the  integrand changes sign.

Consider now a sequence $\breve{K}_m$ of curvatures that converge in $C^3$ to a \emph{forbidden} 
function on $S^3$ or on $S^4$, monotone and non-decreasing in the last Euclidean variable. One could then 
use the analysis in \cite{[ACGPY]} and in \cite{yy2} in dimensions 
three and four respectively to show that blow-ups are isolated and simple near the north pole, 
reaching then a contradiction to existence via the identity \eqref{eq:poho}. 

Applying this reasoning to arbitrarily pinched functions as in \cite{str05} having more than one critical 
point with negative Laplacian, one sees that the dimensional assumption in (ii) of Theorem \ref{t:pinching}
is indeed sharp. 
\end{remark}


\section{Appendix}

Here we collect the proofs of a proposition and  two technical lemmas from the previous sections.

\begin{proof}[Proof of Proposition \ref{p:Km}]
 We illustrate the construction  
dividing it into seven steps. 

\medskip

\noindent {\bf Step 1.} 
Near the south pole $\mathtt{S}$
we can use $\pi_{\mathtt{N}}$ coordinates 
$\{ y_1, \dots, y_n \}$, i.e. coordinates induced by the stereographic projection from the north pole $\mathtt{N}$ mapping 
$\mathtt{S}$ to  $0\in\R^n$. 
For $\d_0 > 0$ and $\e_0 > 0$ small consider a function $\mathcal{K}$ satisfying
$$
\begin{cases}
\mathcal{K} =  \frac{\e_0}{8 n^4} {y}_n^2 
& 
\text{ for } \quad
\mathtt{y}_{n+1} \leq - 1 + \d_0; 
\\ 
\mathcal{K} =  \e_0 (1 + \mathtt{y}_{n+1}) 
& 
\text{ for } \quad
\mathtt{y}_{n+1} \geq - 1 + 2 \d_0; 
\\
\langle \n \mathcal{K}, \n \mathtt{y}_{n+1} \rangle \geq 0 
& 
\text{ on } \quad
S^n \setminus \{ \mathtt{N}, \mathtt{S} \}. 
\end{cases} 
$$
We can also assume that 
$$
\{ \n \mathcal{K} = 0 \} 
\cap 
\{ \mathtt{y}_{n+1} \leq - 1 + 2 \d_0 \} 
\subseteq 
\{ \mathtt{y}_{n+1} \leq - 1 + \d_0 \}
.$$ 
The above function can be chosen so that its Laplacian with respect to the $y$-coordinates 
is bounded away from zero in the set 
$$\{ \mathtt{y}_{n+1} \geq - 1 + 2 \d_0 \}.$$ 
If $\varphi_{\pi_{\mathtt{N}}}$ is the conformal factor of t $\pi_{\mathtt{N}}$ , i.e. 
$g_{S^n} = \varphi_\pi dy^2$,  then 
$$
\Delta_{g_{S^{n}}} \mathcal{K} 
= 
\varphi_\pi^{-1} \Delta_{g_{\R^n}} \mathcal{K} 
+ 
O(\vert \n \varphi_\pi\vert  \, \vert \n \mathcal{K}\vert )
.$$ 
As a consequence $\mathcal{K}$ satisfies 
$$
\D_{g_{S^n}} \mathcal{K} 
\geq 
c > 0
\quad \text{ on } \quad 
U := \{ \mathtt{y}_{n+1} < -1 + 2 \d_0 \}.$$

\medskip

\noindent {\bf Step 2.} We consider next a Morse function $\tilde{K}$ with prescribed numbers of 
critical points with fixed indices and only one local maximum, which we can assume to coincide 
with $\mathtt{N}$. We compose $\tilde{K}$ on the right with a M\"obius map $\Phi$ preserving $\mathtt{N}$ so that 
all other critical points $\{p_1, \dots, p_l\}$ of $\tilde{K} \circ \Phi$ lie in the set 
$\{ \mathtt{y}_{n+1} \leq -1 + \frac{1}{4} \d_0 \}$, 
where $\d_0$ is as in the previous step. The composition with the map $\Phi$ does not affect the 
Morse structure of the function $\tilde{K}$.

\medskip

\noindent {\bf Step 3.} 
For $\d_0$ small the coordinates of the points $p_i$, which we still denote by $p_i$, are of the 
form 
$$
p_i 
= 
(p'_i, p_i^n)
\quad \text{ with } \quad 
p'_i \in \R^{n-1}, p_i^n \in \R
\quad \text{ and } \quad 
(p'_i, p_i^n) \in B_{\d_0^{\frac 14}}(0) \subseteq \R^n. 
$$
By a proper rotation around $0\in\R^n$ we may assume that  
$p'_i \neq  p'_j \in \R^{n-1}$ for $i\neq j$.

\medskip

\noindent {\bf Step 4.} 
Since $\tilde{K} \circ \Phi$ is Morse, 
there exists a rotation $R_i \in SO(n)$ and a diagonal non-singular matrix $A_i$ 
such that near $p_i$ 
$$
(\tilde{K} \circ \Phi)(y) 
= 
\langle R_i (y-p_i), A_i \, R_i (y - p_i) \rangle 
+ 
O(\vert y - p_i\vert ^3).
$$
Without affecting the Morse structure of $\tilde{K}$ we can modify it so that 
one has exactly 
$$
(\tilde{K} \circ \Phi)(y) 
= 
\langle R_i [y-p_i], A_i \, R_i [y - p_i] \rangle  
\quad \text{ for } \quad
\vert y - p_i\vert  \leq \d_1
$$
for some 
$\d_1 \ll \d_0$. Since no $p_i$ is a local maximum, 
we can also assume that the last diagonal entry of $A_i$ is positive.

\medskip

\noindent {\bf Step 5.} 
We next consider a smooth curve $\gamma_{i} : [0,1] \longrightarrow SO(n)$ such that 
$$
\gamma_{i}(0) = Id_{n}
\quad \text{ and } \quad 
\gamma_{i}(1) = R_i, 
$$
and then introduce the new function 
$$
\Theta_i(y) 
: = 
\langle 
\gamma_i (f(\vert y-p_i\vert ^2))  [y-p_i]
,
A_i \, \gamma_i (f(\vert y-p_i\vert ^2))[y - p_i] 
\rangle  
\quad \text{ for } \quad
\vert y - p_i\vert  \leq \d_1, 
$$
where $f$ is zero in a neighbourhood of zero and equal to $1$ in a neighbourhood of $\d_1^2$. We claim that $p_i$
is the only critical point of this function in $B_{\d_1}(p_i)$. In fact consider a curve in $\R^n$ of the type
$$
Y_t 
:= 
p_i 
+  
t (\gamma_i(f(t^2)))^{-1} Y
\quad \text{ with } \quad 
Y \in \R^n, \vert Y\vert  = 1
\quad \text{ and for } \quad 
t \in [0,\d_1].
$$ 
Then clearly 
$\Theta_i(Y_t) =  t^2 \langle Y, A \, Y \rangle$, 
so whenever 
$\langle Y, A \, Y \rangle \neq 0$
the gradient of $\Theta_i$ is non-zero for $t \neq 0$. 
If instead 
$\langle Y, A \, Y \rangle = 0$, 
one can always consider a trajectory $Y_s$ 
in the unit sphere such that 
$$\frac{d}{ds}\lfloor_{s=0} \langle Y_s, A \, Y_s \rangle \neq 0.$$ 
If $Y_t$ is as in the previous formula, consider the curve $Y_t(s)$ replacing $Y$ with $Y(s)$. 
Then its $s$-derivative is a non-critical direction for $\Theta_i$. 
In this way we have proved 
$$
\exists \;
0 < \d_2 \ll \d_1
\;:\;
\Theta_i(y)
= 
\langle 
y-p_i
, 
A_i \,  [y - p_i] 
\rangle  
\quad \text{ for } \quad
\vert y - p_i\vert  \leq \d_2
$$
with diagonal $A_{i}$ and $(A_{i})_{nn}>0$.  
Replacing $\tilde{K} \circ \Phi$ with $\Theta_i$ near each $p_i$, 
no further critical point is created and the Morse structure preserved.


\medskip

\noindent {\bf Step 6.} 
Recall that we rotated the coordinates so that the first $n-1$ components of the 
points $p_i$, i.e. $p'_1, \dots, p'_l \in \R^{n-1}$ are all distinct. 
There exists then 
$$ 
\; \exists \; 0 < \d_3 \ll \d_2
\; \forall \; i\neq j \; : \; \vert p'_i-p'_j\vert  \geq 4 \d_3
$$ 
We choose next a cut-off function $\mathcal{G}$ such that 
$$
\begin{cases}
\mathcal{G} = p_i^n 
&
\text{ in } \quad 
B_{\d_3}(p_i) 
\\ 
\mathcal{G} = 0 
& 
\text{ in } \quad
\R^{n-1} \setminus \cup_{i=1}^l B_{2 \d_3}(p_i).  
\end{cases}  
$$
Calling $\Theta$ the function obtained from replacing 
$\tilde{K}\circ \Phi$ 
by  $\Theta_i$ near $p_i$, we let 
$$
\tilde{\Theta}(y',y_n) 
= 
\Theta(y', y_n + \mathcal{G}(y')). 
$$
Then the only critical points of $\tilde{\Theta}$ are precisely 
$(p'_1, 0), \dots, (p'_l, 0)$. 
In fact these are critical points by construction and moreover 
$$
\begin{cases}
\nabla_{y'} \tilde{\Theta}(y',y_n) 
= 
\nabla_{y'} \Theta(y', y_n + \mathcal{G}(y')) 
- 
\pa_{y_n} \Theta(y', y_n + \mathcal{G}(y')) \nabla_{y'} \mathcal{G}(y'); 
\\  
\pa_{y_n} \tilde{\Theta}(y',y_n) 
= 
\pa_{y_n} \Theta(y', y_n + \mathcal{G}(y')). 
\end{cases}
$$
This implies that
$\nabla \tilde{\Theta}(y',y_n) = 0$ 
if and only if 
$\nabla \Theta(y', y_n + \mathcal{G}(y')) = 0$, 
which is the desired claim.

\medskip
 
\noindent {\bf Final step.}  
Let us call $\hat{K}$ the function obtained from $\tilde{K}$ following the previous steps and consider a sequence of M\"obius maps $\Phi_m$ fixing $\mathtt{N}$ and 
and sending every other point to $\mathtt{S}$ as $m \longrightarrow \infty$. 
Given a Morse function $\tilde{K}$ as in the statement of the proposition, 
we apply the previous steps {\bf{3}}-{\bf{6}}. For $\e_0$ small and fixed and $\e_m \searrow 0$
we then consider a function $K_m$ of the form ($\hat{K} = \hat{K}_m$)
$$
K_m = 1 + \e_0 \mathcal{K} + \e_m \hat{K}_m. 
$$
Using the fact that 
$\mathcal{K} \equiv 0$ for $y_n = 0$ and $\vert y\vert $ small,  
one can check that all critical points of $K_m$ are either at $\mathtt{N}$ 
as the global maximum 
or converge to $\mathtt{S}$  with 
$$
\mathcal{M}_j(K_m) 
= 
\mathcal{M}_j(\tilde{K})
\quad \text{ for all } \quad 
j
$$ 
If $\e_m \searrow 0$ sufficiently fast, 
then $K_m$ satisfies the desired properties with 
$$K_0 = 1 + \e_0 \mathcal{K}.$$ 
\end{proof}

\

\begin{proof}[Proof of Lemma \ref{l:outside}]
We are going to prove the statement using comparison principles on a suitable subset of the sphere. 
First let $G_{\mathtt{N}}$ denote the Green's function of $L_{g_{S^n}}$ with pole at 
$\mathtt{N}$ 
(
$
G_{\mathtt{N}}(x) 
\simeq 
d_{S^n}(x,\mathtt{N})^{2-n}
$ 
near 
$\mathtt{N}$
),  
let $\a \in (0,1)$ and $\d > 0$. 
By direct computation we have that 
\begin{equation}\label{eq:GNa}
\begin{split}
(
L_{g_{S^n}} 
- 
\d d_{S^n}(x,\mathtt{N})^{-2}
) 
(G_{\mathtt{N}})^\a 
= &
\left[ 
(1 - \a) R_{g_{S^n}} 
- 
\d d_{S^n}(x,\mathtt{N})^{-2} 
\right] 
(G_{\mathtt{N}})^\a  \\
& + 
c_n \a (1-\a) (G_{\mathtt{N}})^{\a-2} \vert \n G_{\mathtt{N}}\vert ^2.  
\end{split}
\end{equation}
Fixing first $\a \in (0,1)$ and then $\d > 0$ sufficiently small, 
the right-hand side of \eqref{eq:GNa} is positive.  
Moreover by definition of $\bar{\mu}_m$ and, since $(K_m)_{m}$ is uniformly bounded,
\begin{equation}\label{K_m*u_m_bound}
\exists \; C=C_{\delta}>0
\;:\;
K_m(x) u_m(x)^{\frac{4}{n-2}} 
\leq 
\d d_{S^n}(x,\mathtt{N})^{-2}
\quad \text{ for } \quad  
d_{S^n}(x,\mathtt{N}) \geq C \bar{\mu}_m.
\end{equation}  
In fact, if this inequality were false, from the convergence of $\bar{v}_m$
and the upper bound in \eqref{eq:ubuuu} we could obtain a non-zero limit in 
$C^2_{\text{loc}}(\R^n \setminus \{0\})$ 
for a sequence of the form 
$$\hat{\mu}_m^{\frac{n-2}{2}} u_m(\bar \hat{\mu}_m x)
\quad \text{ with } \quad 
\hat{\mu}_m \ll \bar{\mu}_m,
$$
violating property (ii) before Lemma \ref{l:reg-bubble-limit}. 
Hence \eqref{K_m*u_m_bound} is proved, whence from \eqref{eq:sc-m}
$$
\begin{cases}
(L_{g_{S^n}} - \d d_{S^n}(x,\mathtt{N})^{-2})  u_m 
\leq 
0 
& 
\text{ in } \quad
\{ d_{S^n}(x,\mathtt{N}) \geq C \bar{\mu}_m \}; 
\\ 
u_m \leq \d (C \bar{\mu}_m)^{\frac{2-n}{2}} 
& 
\text{ on } \quad
\{ d_{S^n}(x,\mathtt{N}) = C \bar{\mu}_m \}, 
\end{cases}
$$
while  $G_{\mathtt{N}}$ by \eqref{eq:GNa} is a super-solution of the latter problem on
$\{ d_{S^n}(x,\mathtt{N}) \geq C \bar{\mu}_m \}$. 
By  Hardy-Sobolev's inequality \cite{ckn} and domain monotonicity the quadratic form 
$$
\int_{d_{S^n}(x,\mathtt{N}) \geq C \bar{\mu}_m} v (L_{g_{S^n}} v - \d d_{S^n}(x,\mathtt{N})^{-2} v) d\mu_{g_{S^n}}
$$
is for $\d$ small uniformly positive definite on functions 
vanishing at the boundary of the corresponding spherical cap. 
As a consequence we have a positive first Dirichlet eigenvalue of 
\begin{equation*}
\begin{split}
L_{g_{S^n}} - \d d_{S^n}(x,\mathtt{N})^{-2}
\quad \text{ on } \quad  
\{  d_{S^n}(x,\mathtt{N}) \geq C \bar{\mu}_m\}
\end{split}
\end{equation*}
and this operator satisfies the maximum principle, cf. 
\cite{pw}, \S 5.2, Theorem 10.  
Thus
$$
u_m 
\leq 
(C \bar{\mu}_m)^{\frac{2-n}{2}} 
\left( 
G_{\mathtt{N}}\lfloor_{\pa B_{C \bar{\mu}_m}(\mathtt{N})} 
\right)^{-\a}
G_{\mathtt{N}}^\a  
\leq 
(C \bar{\mu}_m)^{\frac{2-n}{2}} 
\left( 
\frac{\bar{\mu}_m}{d_{S^n}(x,\mathtt{N})} 
\right)^{\a (n-2)}
$$
on
$$
\{ d_{S^n}(x,\mathtt{N}) \geq C \bar{\mu}_m \}. 
$$
Note that $G_{\mathtt{N}}$ is axially symmetric around $\mathtt{N}$, i.e. 
$G_{\mathtt{N}}\lfloor_{\pa B_{C \bar{\mu}_m}(\mathtt{N})}$. 
Hence from \eqref{eq:sc-m}   
\begin{equation}\label{eq:nonlinear}
\begin{cases}
L_{g_{S^n}}u_m 
\leq 
C \bar{\mu}_m^{\frac{2-n}{2}} 
\left(
\frac{\bar{\mu}_m}{d_{S^n}(x,\mathtt{N})} 
\right)^{\a (n+2)}  
&
\text{ in } \quad
\{ d_{S^n}(x,\mathtt{N}) \geq C \bar{\mu}_m \}; 
\\ 
u_m 
\leq 
\d (C \bar{\mu}_m)^{\frac{2-n}{2}} 
& 
\text{ on } \quad
\{ d_{S^n}(x,\mathtt{N}) = C \bar{\mu}_m \}. 
\end{cases}
\end{equation}
%
%
We set 
$
\psi(G_{\mathtt{N}}) 
= 
\Lambda
+ 
\beta (G_{\mathtt{N}})^\gamma
$ 
with 
$\Lambda, \b > 0$ and $\gamma > 1$. 
By direct computation we find
\begin{equation}\label{eq:aaaa}
L_{g_{S^n}} (\psi(G_{\mathtt{N}})) 
= 
\Lambda R_{g_{S^n}} 
+
\beta (\gamma-1) G_{\mathtt{N}}^\gamma 
\left[ 
c_n \gamma \frac{\vert \n G_{\mathtt{N}}\vert ^2}{G_{\mathtt{N}}^2} - R_{g_{S^n}} 
\right]. 
\end{equation}
For $\a < 1$ but close to $1$, we choose $\gamma$ to satisfy 
$$
(n-2) \gamma = \a (n+2) - 2. 
$$ 
Near $\mathtt{N}$ then   
$$
G_{\mathtt{N}}^\gamma 
\frac
{\vert \n G_{\mathtt{N}}\vert ^2}
{G_{\mathtt{N}}^2}
\sim
d_{S^n}(x,\mathtt{N})^{- \a (n+2)},$$ 
as is the right-hand side of the first inequality in \eqref{eq:nonlinear}, while 
$G_{\mathtt{N}}^\gamma R_{g^{S^n}}$ 
is of lower order. Choosing $\beta$ to satisfy 
$$
\beta \bar{\mu}_m^{-\a (n+2)} 
= 
\bar{C} \bar{\mu}_m^{\frac{2-n}{2}}
\quad \text{ with } \quad 
\bar{C} \gg C 
\quad \text{ large fixed}, 
$$
near $\mathtt{N}$ the right-hand side in \eqref{eq:aaaa} dominates the one in \eqref{eq:nonlinear}. 
Choosing in addition 
$$
\beta \ll \Lambda\ll \mu^{\frac{n-2}{2}}, 
$$
which is possible by the above choice of $\b$, then we obtain the properties 
$$
\begin{cases} 
\Lambda
+ 
\beta (G_{\mathtt{N}})^{\gamma} 
\leq 
C \bar{\mu}_m^{\frac{n-2}{2}} d_{S^n}(x,\mathtt{N})^{2-n} 
& 
\text{ in } \quad
\{ d_{S^n}(x,\mathtt{N}) \geq C \bar{\mu}_m \}; 
\\ 
L_{g_{S^n}} (\Lambda + \beta (G_{\mathtt{N}})^{\gamma}) 
\geq 
L_{g_{S^n}} u_m 
& 
\text{ in } \quad
\{ d_{S^n}(x,\mathtt{N}) \geq C \bar{\mu}_m \}; 
\\ 
u_m 
\leq 
\Lambda 
+ 
\beta (G_{\mathtt{N}})^{\gamma} 
& 
\text{ on } \quad
\{ d_{S^n}(x,\mathtt{N}) = C \bar{\mu}_m \}. 
\end{cases}
$$
Then the conclusion follows from the maximum principle. 
\end{proof}
 
\begin{proof}[Proof of Lemma \ref{l:sing-prof-anyway}]
We follow the proof of Proposition 2.3 in \cite{yy}, which relies on 
Proposition 2.1, Lemma 2.1, Lemma 2.3 and Lemma 2.3 there.
The crucial point here is that  uniform 
gradient bounds on $\check{K}_m$ fail, so we cannot directly extract a bubble from the maximum 
point of $\check{u}_m$. We can however exploit the estimate in Lemma \ref{l:outside} instead. 
Apart from some modifications that we will describe in detail, the arguments there can 
be carried out even without gradient bounds.

Similarly to \cite{yy} consider a maximum point $y_m$ of $\check{u}_m$, 
a unit vector $e \in \R^n$ and 
$$
\check{v}_m(y) 
= 
\check{u}_m(y_m + e)^{-1} \check{u}_m(y)
.
$$ 
As in there we prove that 
$\check{v}_m$ converges in $C^2_{\text{loc}}(B_1\setminus \{0\})$
to a singular function 
$$\check{v}(y) = a \vert y\vert ^{2-n} + h(y)$$ 
with $a > 0$ and $h$ smooth and harmonic. 
The next step consists in showing that 
\begin{equation}\label{eq:yy_to_be_shown}
\check{u}_m(y_m + e) 
\leq 
C \check{u}_m(y_m)^{-1} 
\end{equation}
for some fixed $C > 0$. If this is not true, then we have 
\begin{equation}\label{eq:5.4}
\limsup_m \check{u}_m(y_m) \check{u}_m(y_m + e) 
\longrightarrow 
\infty. 
\end{equation}
Multiplying \eqref{eq:sat} by 
$\check{u}_m(y_m+e)^{-1}$  one finds after integration
\begin{equation*}
\begin{split}
- \oint_{\pa B_1} \frac{\pa}{\pa \nu} \check{v}_m d \s 
= &
- \check{u}_m(y_m+e)^{-1} \int_{B_1} \Delta \check{u}_m dx \\
= &
\frac{1}{c_n} \check{u}_m(y_m + e)^{-1} 
\int_{B_1} \check{K}_m \check{u}_m^{\frac{n+2}{n-2}} dx. 
\end{split}
\end{equation*} 
From the fact that $h$ is harmonic and that $a > 0$ we get that 
$$
\lim_m \oint_{\pa B_1} \frac{\pa}{\pa \nu} \check{v}_m d \s 
= 
\oint_{\pa B_1} 
\frac{\pa}{\pa \nu} \left( a \vert y\vert ^{2-n} + h(y) \right) d \s 
< 
0. 
$$
For
$R_m \longrightarrow \infty$ 
sufficiently slowly set 
$$r_m = R_m \check{u}_m(y_m)^{-\frac{2}{n-2}}.$$ 
Then by Lemma \ref{l:outside} and a change of variables  
$$
\int_{\vert y - y_m\vert  \leq r_m} \check{K}_m \check{u}_m^{\frac{n+2}{n-2}} dx 
\leq 
C \check{u}_m(y_m)^{-1}. 
$$
As for Lemma 2.2 in \cite{yy}, which is based on local estimates in the annulus 
$$r_m \leq \vert  y - y_m \vert  \leq 1$$ only,  
it is possible to prove that 
$$
\check{u}_m(y) 
\leq 
C \check{u}_m(y_m)^{-\check{\l}_m} \vert y - y_m\vert ^{2-n+\d_m} 
\quad \text{ for } \quad
r_m \leq \vert  y - y_m \vert  \leq 1, 
$$ 
where 
$\d_m = O(R_m^{-2+o_m(1)})$
and 
$\check{\l}_m = \frac{2(n-2-\d_m)}{n-2}-1$. 
This implies 
$$
\int_{r_m \leq \vert y- y_m\vert  \leq 1} \check{K}_m \check{u}_m^{\frac{n+2}{n-2}} dx 
\leq 
C R_m^{n-\frac{n+2}{n-2} (n-2-\d_m)} \check{u}_m(y_m)^{-1} 
= 
o(1) \check{u}_m(y_m)^{-1}.  
$$
The latter formulas would then give a contradiction to \eqref{eq:5.4}.
Hence 
\eqref{eq:yy_to_be_shown} is established and the rest of the proof of Proposition 2.3 in \cite{yy} goes through in our case too. 
\end{proof}

\

\ 

\hline
\begin{center}
Acknowledgments
\end{center}	
\begin{enumerate}[label=(\roman*)]
\item 	A.Malchoidi has been supported by the project {\em Geometric Variational Problems} and {\em Finanziamento a supporto della ricerca di base} from Scuola Normale Superiore and by MIUR Bando PRIN 2015 2015KB9WPT$_{001}$.  He is also member of GNAMPA as part of INdAM. 
\item M.Mayer has been supported by the Italian MIUR Department of Excellence grant CUP E83C18000100006. 
\end{enumerate}

\bibliographystyle{cpam}



                                

\end{document}